\documentclass[a4paper,leqno]{amsart}

\usepackage{amsmath,amssymb,amsthm}
\usepackage{hyperref}

\def\H{{\mathbb H}}% hyperbolic space
\def\C{{\mathbb C}}% complex numbers
\def\R{{\mathbb R}}% real numbers
\def\Sph{{\mathbb S}} % sphere
\def\F{\mathcal F}% Fourier transform
\def\M{\mathcal M}

\def\coth{\operatorname{cotanh}}
\def\virgp{\raise 2pt\hbox{,}}

\def\({\left(}
\def\){\right)}
\def\<{\left\langle}
\def\>{\right\rangle}
\def\le{\leqslant}
\def\ge{\geqslant}

\def\Eq#1#2{\mathop{\sim}\limits_{#1\rightarrow#2}}
\def\Tend#1#2{\mathop{\longrightarrow}\limits_{#1\rightarrow#2}}

\def\d{{\partial}}

\def\eps{\varepsilon}
\def\l{\lambda}
\def\om{\omega}
\def\si{{\sigma}}

\def\w{{\tt w}_n}

\def\gd{{S_d}}

\DeclareMathOperator{\RE}{Re}
\DeclareMathOperator{\IM}{Im}

\theoremstyle{plain}
\newtheorem{theo}{Theorem}[section]
\newtheorem{lem}[theo]{Lemma}
\newtheorem{cor}[theo]{Corollary}
\newtheorem{prop}[theo]{Proposition}

\theoremstyle{definition}
\newtheorem{defin}[theo]{Definition}
\newtheorem*{nota}{Notation}

\theoremstyle{remark}
\newtheorem{rema}[theo]{Remark}
\newtheorem*{rema*}{Remark}

\numberwithin{equation}{section}

\begin{document}

\title[Scattering for NLS on hyperbolic space]{Scattering theory
  for radial nonlinear Schr\"odinger equations on hyperbolic space}  
\author[V. Banica]{Valeria Banica}
\address[V. Banica]{D\'epartement de Math\'ematiques\\ Universit\'e
  d'Evry\\ Bd. F.~Mitterrand\\ 91025 Evry\\ France} 
\email{Valeria.Banica@univ-evry.fr}
\author[R. Carles]{R{\'e}mi Carles}
\address[R. Carles]{MAB, UMR 5466 CNRS\\Universit\'e de Bordeaux 1\\
  351 cours de la  
  Lib{\'e}ration\\ 33405 Talence cedex\\ France}
\email{Remi.Carles@math.cnrs.fr}
\author[G. Staffilani]{Gigliola Staffilani}
\address[G. Staffilani]{MIT\\77 Massachusetts Avenue\\ Cambridge\\ MA
  02139-4307\\ USA} 
\email{gigliola@math.mit.edu}
\thanks{V.B. is partially
  supported by the ANR project ``\'Etude qualitative des E.D.P.''.
R.C. acknowledges supports by European network HYKE,
  funded  by the EC as contract HPRN-CT-2002-00282, and by the ANR
  project SCASEN. G. S. is partially supported by N.S.F. Grant  0602678.} 
\begin{abstract}
We study the long time behavior of radial solutions to nonlinear
Schr\"odinger equations on hyperbolic space. We show that the
usual distinction between short range and long range nonlinearity is
modified: the geometry of  the hyperbolic space makes every power-like
nonlinearity short range. The proofs rely on weighted Strichartz
estimates, which imply Strichartz estimates for a broader
family of admissible pairs, and on Morawetz type inequalities. The
latter are established without symmetry assumptions. 
\end{abstract}
\subjclass[2000]{}
\maketitle
\tableofcontents

\section{Introduction}
\label{sec:intro}
This paper is devoted to the scattering theory for the nonlinear
Schr\"odinger equation
\begin{equation}
  \label{eq:nls}
i\d_t u + \Delta_{\H^n} u=|u|^{2\si}u \quad ;\quad
  e^{-it\Delta_{\H^n}}u(t)\big|_{t=t_0}=\varphi\, ,
\end{equation}
on the hyperbolic space ($n\ge 2$):
\begin{equation*}
  %\label{eq:H}
  \H^n = \{ \R^{n+1}\ni \Omega =(x_0,\ldots,x_n)=(x_0,x')=(\cosh
  r,\sinh r\, 
  \om),\ r\ge 0, \ \om \in \Sph^{n-1} \}.
\end{equation*}
We consider  a defocusing power nonlinearity.
One  could also prove some results in the focusing case,
but this case will not be  discussed in this paper. 
When a function of time and space $u(t,\Omega)$ depends only on $t$
and $r$, we say that it is radial. The reason is that $r$ is the
hyperbolic distance between $\Omega$ and the origin of the hyperboloid
$O=(x_0=1,x'=0)$. With the usual abuse of notation, we write
$u(t,r)$. We prove that for any $\si>0$, a short range (large data)
scattering theory is 
available for radial solutions to \eqref{eq:nls}. This is in sharp
contrast with the Euclidean case, where the nonlinearity $|u|^{2\si}u$
cannot be short range as soon as 
$\si\le 1/n$ (see Section~\ref{sec:longrange}). A crucial argument to
prove this phenomenon is the existence of weighted Strichartz
estimates for radial solutions to Schr\"odinger equations on $\H^n$,
established in \cite{BaHyper} and \cite{VittoriaDR}. Note that if
these weighted Strichartz
estimates were available for general solutions to Schr\"odinger
equations on $\H^n$ (and not only radial), then all the results of
this paper could be adapted, with the same proofs. Also, similar
results can be extended to the equation posed on Damek-Ricci spaces,
thanks to the weighted Strichartz estimates obtained in
\cite{VittoriaDR}.  
Finally, let us recall that recently the nonlinear Schr\"odinger equation in
a non-Euclidean setting has been intensively studied
(see e.g. \cite{BourgainTorus,BGTMRL,BGT}). 
Most of the results concern the local-in-time point of view, and to
the best of our knowledge, until now there was no result of large data
scattering in a non-Euclidean manifold.  
\smallbreak
With the above parameterization for  the hyperbolic space, the
Laplace-Beltrami operator reads:
\begin{equation*}
  %\label{eq:Laplace}
  \Delta_{\H^n}= \d_r^2 + (n-1)\frac{\cosh r}{\sinh r}\, \d_r +\frac{1}{\sinh^2
  r}\Delta_{\Sph^{n-1}}\, .
\end{equation*}
In order to define wave operators,  we introduce the free
Schr\"odinger generalized initial value problem:  
\begin{equation}
  \label{eq:linear}
\left\{\begin{array}{c}  i\d_t u + \Delta_{\H^n} u=0\quad,\\
\quad u_{\mid t=t_0}=\, \varphi.\end{array}\right.
\end{equation}
We denote $U(t)= e^{it\Delta_{\H^n}}$, so that in
\eqref{eq:linear}, $u(t,\Omega) =U(t)\varphi(\Omega)$.
When considering solutions to \eqref{eq:nls}, we use the convention that
if $t_0=-\infty$ (resp. $t_0=+\infty$), 
then we denote 
$\varphi=u_-$ (resp. $\varphi=u_+$), and solving
\eqref{eq:nls}  means that we construct wave
operators. If $t_0=0$, then we denote $\varphi =u_0$, and
\eqref{eq:nls} is the standard Cauchy
problem. In all the cases, we seek mild solutions to
\eqref{eq:nls}, that is, we solve
\begin{equation}
  \label{eq:duhamel}
  u(t)= U(t)\varphi -i\int_{t_0}^t U(t) (t-s)\(|u|^{2\si}u\)(s)ds.  
\end{equation}
We can now state our main  results. The first one deals with existence
of wave operators and asymptotic completeness for small 
$L^2$ data:
\begin{theo}\label{theo:smallL2}
  Let $n\ge 2$, $0<\si\le 2/n$, and $t_0\in \overline\R$. There exists
  $\epsilon=\epsilon(n,\si)$ such that if $\varphi\in L^2_{\rm rad}(\H^n)$ 
  with $\|\varphi\|_{L^2}<\epsilon$, then 
  \eqref{eq:nls}  has a unique solution
  \begin{equation*}
    u\in C(\R;L^2)\cap L^{2+2\si}(\R\times\H^n).
  \end{equation*}
Moreover, its $L^2$-norm is constant, $\|u(t)\|_{L^2}=
\|\varphi\|_{L^2}$ for all $t\in \R$.\\ 
There exist $u_\pm\in L^2_{\rm
  rad}(\H^n)$ 
such that
\begin{equation*}
  \left\| u(t)-U(t)u_\pm\right\|_{L^2}\to 0 \quad \text{as }t\to \pm
  \infty. 
\end{equation*}
If $t_0=-\infty$ (resp. $t_0=+\infty$), then $u_-=\varphi$
(resp. $u_+=\varphi$). 
\end{theo}
The existence of solutions in $C(\R;L^2)$ for data which are small in
$L^2$ is analogous to  the Euclidean case
(\cite{TsutsumiL2}, see also \cite{CazCourant}). For $\si=2/n$, our
result is the exact analogue to its Euclidean counterpart recalled in
Proposition 2.3. Note 
however that for $0<\si<2/n$, the space where the solutions belong,
and the existence of a scattering theory, distinguish the hyperbolic
space $\H^n$ from the Euclidean space $\R^n$. In particular, there is
no long range effect in hyperbolic space, even if $0<\si\le 1/n$.
\smallbreak

Our second result establishes the existence of  the wave operator in the
Sobolev space $H^1$, when the nonlinearity is $H^1$-subcritical (see
Appendix~\ref{sec:crit} for the notion of criticality). Here again,
the power $\si$ can go down to $0$, with no long range effect. 
\begin{theo}\label{theo:waveH1}
   Let $n\ge 2$, $0<\si<2/(n-2)$, and $t_0=-\infty$. For any $\varphi=u_-\in
  H^1_{\rm rad}(\H^n)$, 
  there exists $T<\infty$ such that \eqref{eq:duhamel} has a
  unique solution in $C\cap L^\infty(]-\infty,-T];H^1)\cap
  L^{2\si+2}(]-\infty,-T];W^{1,2\si +2}) $. \\
Moreover, this solution $u$ is defined globally in time:  $u\in
L^\infty(\R;H^1)$. That is, $u$ is the only solution to \eqref{eq:nls}
  with
  \begin{equation*}
    \left\|u(t)-U(t)u_-\right\|_{H^1}=
    \left\|U(-t)u(t)-u_-\right\|_{H^1}\to 0 \quad \text{as }t\to
    -\infty. 
  \end{equation*}
\end{theo}
Of course, we could prove the existence of wave operators with data at
time $t_0=+\infty$. Since the proof is similar, we shall skip it. 
\smallbreak

The proofs of Theorems~\ref{theo:smallL2} and \ref{theo:waveH1} rely on
two remarks. First, the weighted Strichartz estimates proven in
\cite{BaHyper,VittoriaDR} for radial solutions to Schr\"odinger
equations on $\H^n$, $n\ge 3$, make it possible to state 
Strichartz estimates which are the same as on $\R^d$, for any $d\ge
n$. We show in this paper that similar results are available when
$n=2$. Second, the classical proofs for the counterparts of
Theorems~\ref{theo:smallL2} and \ref{theo:waveH1}  in the Euclidean
space $\R^d$ 
rely only on functional
analysis arguments, based on Strichartz estimates, H\"older inequality
and Sobolev embeddings. This is why the proofs of
Theorems~\ref{theo:smallL2} and \ref{theo:waveH1}, presented in
Sections~\ref{sec:small} and \ref{sec:wave} respectively, are rather
short. Finally, let us notice that in dimension $3$, the Strichartz
estimates without weights were proved to hold also for non-radial data
(\cite{BaHyper}). Therefore Theorems~\ref{theo:waveH1} holds for the
usual range of nonlinearities $2/3<\sigma<2$ without symmetry
assumption. 
\smallbreak

The next natural step in scattering theory consists in proving the
invertibility of the wave operators on their range, that is,
asymptotic completeness. We prove a Morawetz type inequality that
combined with the Strichartz estimates for higher-dimension admissible
couples gives us a scattering result without lower restriction on the
nonlinearity power. 
%Classically for nonlinear Schr\"odinger
%equations, this relies on Morawetz type inequalities (\cite{GV85}, see
%also \cite{CazCourant}). However, the original proof uses dispersive
%inequalities, and we could not find a better time decay in $\H^n$ than
%in $\R^n$. This is why we use the more recent approach of
%\cite{TaoVisanZhang}, based on an interaction Morawetz
%inequality. This interaction inequality for hyperbolic space was
%noticed in \cite{HTW05}. 
Note that the asymptotic completeness
that we prove is for
$n=3$ only (see the discussion in Section~\ref{sec:completeness}). 
\begin{theo}\label{theo:CA}
  Let $n=3$, $0<\si<2$, and $t_0=0$. For any $\varphi=u_0\in H^1_{\rm
  rad}(\H^3)$  
  \eqref{eq:nls} has a unique, global solution in $C(\R;H^1_{\rm
  rad})\cap L^4(\R\times \H^3)$. Moreover, there exist $u_-$ and $u_+$
  in $H^1_{\rm rad}(\H^3)$ such that
  \begin{equation*}
    \left\| u(t)- U(t)u_\pm\right\|_{H^1(\H^3)}\to 0 \quad\text{as }
    t\to {\pm \infty}. 
  \end{equation*}
Moreover, if $2/3<\si<2$, then we need not assume that $\varphi$ is
radial: for any $\varphi=u_0\in H^1(\H^3)$,
  \eqref{eq:nls} has a unique, global solution in $C(\R;H^1)\cap
  L^4(\R\times \H^3)$, and there exist $u_-$ and $u_+$ 
  in $H^1(\H^3)$ such that
  \begin{equation*}
    \left\| u(t)- U(t)u_\pm\right\|_{H^1(\H^3)}\to 0 \quad\text{as }
    t\to {\pm \infty}. 
  \end{equation*}
\end{theo}

\begin{nota}
  In this  paper we often use the notation $A\lesssim B$
to denote that there exists an absolute constant $C$ such that $A\le 
C B$.  Another standard notation is to use for any  $1\le p\le
\infty$  the symbol $p'$  to denote the  H\"older-conjugate exponent, that is
$1/p+1/p'=1$. 
%Finally we use the symbol $\RE z$ and $\IM z$ to denote
%respectively the  the real and the imaginary part  of a complex number
%$z$. 
\end{nota}

The rest of this paper is organized as follows. In
Section~\ref{sec:euclidien}, we review the scattering result for
nonlinear Schr\"odinger equations on the Euclidean space: small $L^2$
data, existence of wave operators in $H^1(\R^d)$, non-existence of
wave operators when $\si\le 1/d$, and asymptotic completeness. In
Section~\ref{sec:strichartz}, we show that 
for radial solutions to \eqref{eq:nls}, the same Strichartz estimates
as in $\R^d$ are available in $\H^n$, for any $d\ge n\ge
2$. Theorems~\ref{theo:smallL2} and \ref{theo:waveH1} are proven in 
Sections~\ref{sec:small} and \ref{sec:wave} respectively. We prove a
general interaction Morawetz inequality in Section~\ref{sec:morawetz},
and infer Theorem~\ref{theo:CA} in Section~\ref{sec:completeness}. In
Appendix~\ref{sec:crit} , we prove that the notion of criticality, as far as
the Cauchy problem \eqref{eq:nls} is concerned, is the same on $\H^n$
as on $\R^n$. We study the large time
behavior of 
radial solutions to the linear Schr\"odinger equation
\eqref{eq:linear} on $\H^3$ in Appendix~\ref{sec:DA}. Finally, we
discuss the existence of an analogue to the Galilean operator in the
radial framework on $\H^3$ in
Appendix~\ref{sec:galileo}.

\section{A review of scattering theory in $\R^d$}
\label{sec:euclidien}

In this paragraph, we consider, in the Euclidean space, the equation
\begin{equation}
  \label{eq:nlseucl}
  i\d_t u+\Delta_{\R^d}u = |u|^{2\si}u\quad ;\quad (t,x)\in
  \R\times \R^d.
\end{equation}
We recall some results concerning scattering theory, in order to
compare them with their counterpart in hyperbolic space. We also
sketch some proofs that we mimic in the hyperbolic setting.
\smallbreak

First, the Schr\"odinger operator in the Euclidean space satisfies the
following Strichartz estimates.
\begin{defin}\label{def:adm}
 Let $d\ge 2$. A pair $(p,q)$ is {\em $d$-admissible} if $2\le q
  \le \frac{2d}{d-2}$ and 
  \begin{equation*}
    \frac{2}{p}=\delta(q):= d\left(
    \frac{1}{2}-\frac{1}{q}\right),\quad (p,q)\not = (2,\infty)\, .
  \end{equation*}
\end{defin}
\begin{prop}\label{prop:strieucl}
  Let $d\ge 2$. Denote $\gd(t) = e^{it\Delta_{\R^d}}$. \\
 {\rm1.} For any $d$-admissible pair $(p,q)$, there exists $C_q$
    such that
\begin{equation*}
    \left\| \gd(\cdot)\phi \right\|_{L^p({\R};L^q)}\le C_q \|\phi
    \|_{L^2},\quad \forall \phi \in L^2(\R^d).
  \end{equation*}

\noindent{\rm2.} For any $d$-admissible pairs $(p_1,q_1)$ and $
    (p_2,q_2)$ and any 
    interval $I$, there exists $C_{q_1,q_2}$ independent of $I$ such that 
\begin{equation*}
      \left\| \int_{I\cap\{s\le
      t\}} \gd(t-s)F(s)ds 
      \right\|_{L^{p_1}(I;L^{q_1})}\le C_{q_1,q_2} \left\|
      F\right\|_{L^{p'_2}\(I;L^{q'_2}\)},
    \end{equation*}
for every $F\in L^{p'_2}\(I;L^{q'_2}(\R^d)\)$.
\end{prop}
Let $t_0\in \overline\R$, and consider \eqref{eq:nlseucl} along with
the initial data:
\begin{equation}
  \label{eq:cieucl}
  \gd (-t)u(t)\big|_{t=t_0}=\varphi.
\end{equation}
We use the convention that if $t_0=-\infty$ (resp. $t_0=+\infty$),
then we denote 
$\varphi=u_-$ (resp. $\varphi=u_+$), and solving
\eqref{eq:nlseucl}--\eqref{eq:cieucl} means that we construct wave
operators. If $t_0=0$, then we denote $\varphi =u_0$, and
\eqref{eq:nlseucl}--\eqref{eq:cieucl} is the standard Cauchy
problem. In all these cases, we seek mild solutions to
\eqref{eq:nlseucl}--\eqref{eq:cieucl}, that is, we solve
\begin{equation}
  \label{eq:duhameleucl}
  u(t)= \gd (t)\varphi -i\int_{t_0}^t \gd (t-s)\(|u|^{2\si}u\)(s)ds=:
  \Phi(u)(t).  
\end{equation}
\subsection{Small data in the $L^2$-critical case}
\label{sec:smalleucl}

Recall the result of \cite{CW89}. 
The $L^2$-critical case corresponds to the power $\si =2/d$. In that
case, the pair
\begin{equation*}
  (p,q)= \( 2+\frac{4}{d},2+\frac{4}{d}\).
\end{equation*}
is $d$-admissible, and this is the main remark to prove:
\begin{prop}\label{prop:L2eucl}
  Let $d\ge 2$, $\si=2/d$, and $t_0\in \overline\R$. There exists $\epsilon=\epsilon(d)$ such that if
  $\varphi\in L^2(\R^d)$ 
  with $\|\varphi\|_{L^2}<\epsilon$, then 
  \eqref{eq:nlseucl}--\eqref{eq:cieucl}  has a unique solution
  \begin{equation*}
    u\in C(\R;L^2)\cap L^{2+\frac{4}{d}}(\R\times\R^d).
  \end{equation*}
Moreover, its $L^2$-norm is constant, $\|u(t)\|_{L^2}=
\|\varphi\|_{L^2}$ for all $t\in \R$.\\ There exist $u_\pm\in L^2(\R^d)$
such that
\begin{equation*}
  \left\| u(t)-\gd(t)u_\pm\right\|_{L^2}\to 0 \quad \text{as }t\to \pm
  \infty. 
\end{equation*}
If $t_0=-\infty$ (resp. $t_0=+\infty$), then $u_-=\varphi$
(resp. $u_+=\varphi$). 
\end{prop}
\begin{proof}[Sketch of the proof]
The idea is to apply a fixed point argument to
\eqref{eq:duhameleucl} in 
\begin{equation*}
  X= \left\{u\in C(\R;L^2)\cap L^{2+\frac{4}{d}}(\R\times\R^d)\quad
  ;\quad \|u\|_{L^{2+\frac{4}{d}}(\R\times\R^d)}\le
  2C_{2+\frac{4}{d}}\|\varphi\|_{L^2} \right\}. 
\end{equation*}
Here, $C_{2+\frac{4}{d}}$ is the constant given in the first part of
Proposition~\ref{prop:strieucl}. Indeed, denoting $\gamma=2+4/d$,
Strichartz estimates and H\"older inequality yield:
\begin{equation*}
  \|\Phi(u)\|_{L^{\gamma}(\R\times\R^d)} \le
  C_{\gamma}\|\varphi\|_{L^2(\R^d)} +
  C_{\gamma,\gamma}\|u\|_{L^{\gamma}(\R\times\R^d)}^{1+\frac{4}{d}}. 
\end{equation*}
This shows that for $\|\varphi\|_{L^2}$ sufficiently small, $X$ is
invariant under the action of $\Phi$. Similarly, $\Phi$ is a contraction
on $X$ if $\|\varphi\|_{L^2}$ is sufficiently small, thus providing a
unique solution to \eqref{eq:duhameleucl} in $X$. The conservation of
mass is classical, and holds without the smallness assumption. 
\smallbreak
 
Scattering then
follows from the Cauchy criterion: for $t_1\le t_2$, we have
\begin{equation*}
  \left\| \gd(-t_2)u(t_2)-\gd(-t_1)u(t_1)\right\|_{L^2}\le
  C_{2,\gamma}\|u\|_{L^\gamma ([t_1,t_2]\times \R^d)}^{1+\frac{4}{d}}.
\end{equation*}
The right hand side goes to zero when $t_1,t_2\to \pm\infty$. The
proposition follows easily, since the group $\gd$ is unitary on $L^2$. 
\end{proof}

\subsection{Existence of wave operators in $H^1$}
\label{sec:waveeucl}

We recall the existence of wave operators for negative time; for
positive time, the proof is similar. This means that we solve
\eqref{eq:duhameleucl} with $t_0=-\infty$ (and $\varphi=u_-$). The
strategy consists first in solving \eqref{eq:duhameleucl} in a neighborhood
of $t=-\infty$, that is on $]-\infty,-T]$ for $T$ possibly very
large. Then the conservation of mass and energy makes it possible to
extend the solution to $t\in\R$. We simply recall the first step. The
proof of this result appears in \cite{GV85}. The proof we give is a
simplification, which may be found for instance in \cite{Ginibre}. We
shall not recall or use the results available in weighted Sobolev
spaces (see e.g. \cite{CW92,GOV,GV79Scatt}). 
\begin{prop}\label{prop:waveH1eucl}
  Let $t_0=-\infty$, $d\ge 2$ and $2/d\le\si<2/(d-2)$. For any $\varphi=u_-\in
  H^1(\R^d)$, 
  there exists $T<\infty$ such that \eqref{eq:duhameleucl} has a
  unique solution in $C\cap L^\infty(]-\infty,-T];H^1)\cap
  L^p(]-\infty,-T];W^{1,2\si +2}) $, where $p$ is such that
  $(p,2\si+2)$ is $d$-admissible. \\
Moreover, this solution $u$ is defined globally in time:  $u\in
L^\infty(\R;H^1)$. 
\end{prop}
In other words, we construct the only solution $u$ to
\eqref{eq:nlseucl} such that
\begin{equation*}
  \left\| u(t)-\gd(t)u_-\right\|_{H^1}= \left\|\gd(-t)
  u(t)-u_-\right\|_{H^1}\to 0\quad\text{as }t\to -\infty.
\end{equation*}
The wave operator $W_-$ is the map
\begin{equation*}
  W_-:\quad H^1\ni u_-\mapsto u_{\mid t=0}\in H^1.
\end{equation*}
\begin{proof}
  Recall that $p$ is such that $(p,2\si+2)$ is
  $d$-admissible:
  \begin{equation*}
   p=\frac{4\si+4}{d\si}\cdot 
  \end{equation*}
With the notation $L^\beta_TY=L^\beta(]-\infty,-T];Y)$, we
  introduce:  
  \begin{align*}
    X_T:=\Big\{ u\in C(]-\infty,-T];H^1)\ ;\ &\left\|
    u\right\|_{L^p_TW^{1,2\si+2}} \le 2 C_{2\si+2}\|u_-\|_{H^1},\\
 \left\|  u\right\|_{L^\infty_TH^1} \le 2 \|u_-\|_{H^1}\, ,\quad
& \left\|  u\right\|_{L^p_T L^{2\si+2}} \le 2 \left\|
    \gd(\cdot)u_-\right\|_{L^p_T L^{2\si+2}} \Big\},
  \end{align*}
where $C_{2\si+2}$ is given by Proposition~\ref{prop:strieucl}. Set
$q=s=2\si +2$: we have
\begin{align*}
  \frac{1}{q'}&= \frac{1}{q}+\frac{2\si}{s}\virgp\\
\frac{1}{p'}&= \frac{1}{p}+\frac{2\si}{k}\virgp
\end{align*}
where $(p,q)$ is $d$-admissible and $p\le k<\infty$ since $2/d\le
\si<2/(d-2)$. For $u\in X_T$, Strichartz estimates and H\"older
inequality yield:
\begin{align*}
  \left\| \Phi(u)\right\|_{L^p_T W^{1,2\si +2}} &\le C_{2\si
  +2}\|u_-\|_{H^1} + C\( \left\| |u|^{2\si}u\right\|_{L^{p'}_TL^{q'}}
  + \left\| |u|^{2\si}\nabla u\right\|_{L^{p'}_TL^{q'}}\)\\
&\le C_{2\si
  +2}\|u_-\|_{H^1} + C\|u\|_{L^k_TL^s}^{2\si}\( \|u\|_{L^p_T L^q}
  +\|\nabla u\|_{L^p_T L^q} \)\\
&\le C_{2\si
  +2}\|u_-\|_{H^1} + C\|u\|_{L^p_TL^q}^{2\si\theta
  }\|u\|_{L^\infty_TL^q}^{2\si(1-\theta) } \|u\|_{L^p_T W^{1,2\si
  +2}}\ ,
\end{align*}
for some $0<\theta\le 1$, where we have used the property
$q=s=2\si+2$. Sobolev embedding and the definition of $X_T$ then imply:
\begin{align*}
 \left\| \Phi(u)\right\|_{L^p_T W^{1,2\si +2}} \le C_{2\si
  +2}\|u_-\|_{H^1} + C\left\|\gd(\cdot)u_-\right\|_{L^p_TL^q}^{2\si\theta
  }\|u\|_{L^\infty_TH^1}^{2\si(1-\theta) } \|u\|_{L^p_T W^{1,2\si
  +2}} .
\end{align*}
We have similarly
\begin{align*}
  \left\| \Phi(u)\right\|_{L^\infty_T H^1}&\le  \|u_-\|_{H^1} +
  C\left\|\gd(\cdot)u_-\right\|_{L^p_TL^q}^{2\si\theta 
  }\|u\|_{L^\infty_TH^1}^{2\si(1-\theta) } \|u\|_{L^p_T W^{1,2\si
  +2}} \\
\left\| \Phi(u)\right\|_{L^p_T L^{2\si +2}}& \le
  \left\|\gd(\cdot)u_-\right\|_{L^p_TL^{2\si+2}} 
+ C\left\|\gd(\cdot)u_-\right\|_{L^p_TL^q}^{2\si\theta
  }\|u\|_{L^\infty_TH^1}^{2\si(1-\theta) } \|u\|_{L^p_T W^{1,2\si
  +2}}.
\end{align*}
>From Strichartz estimates, $\gd(\cdot)u_- \in L^p(\R;L^q)$, so 
\begin{equation*}
  \left\|\gd(\cdot)u_-\right\|_{L^p_TL^q} \to 0\quad \text{as }T\to +\infty.
\end{equation*}
Since $\theta>0$, we infer that $\Phi$ sends $X_T$ to itself, for
$T$ sufficiently large. 
\smallbreak

We have also, for $u_2,u_1\in X_T$:
\begin{align*}
  \left\| \Phi(u_2)-\Phi(u_1)\right\|_{L^p_T L^q}&\lesssim
  \max_{j=1,2}\| u_j\|_{L^k_TL^s}^{2\si} \left\|
  u_2-u_1\right\|_{L^p_T L^q}\\
&\lesssim \left\|\gd(\cdot)u_-\right\|_{L^p_TL^q}^{2\si\theta
  }\|u_-\|_{H^1}^{2\si(1-\theta) }\left\|
  u_2-u_1\right\|_{L^p_T L^q}.
\end{align*}
Up to choosing $T$ larger, $\Phi$ is a contraction on $X_T$, and
the proposition follows.
\end{proof}

\subsection{Non-existence of wave operators for $\si\le 1/d$}
\label{sec:longrange}

Even though the scattering result we recalled shows the existence of
wave operators in $H^1(\R^d)$ for $\si \ge 2/d$, it is natural to
expect the nonlinearity to be negligible for large time as soon as
$\si >1/d$. Many results exist, supporting this assertion; we shall
not state them, but rather point out that it is not possible
to go below $1/d$. The result recalled below was established in
\cite{Strauss74,Barab} (see also \cite{HT87a}). 
\begin{prop}\label{prop:longrange}
  Let $d\ge 2$, $0<\si\le 1/d$ and $T>0$. Let $u\in
  C(]-\infty,-T];L^2(\R^d))$ be a solution of \eqref{eq:nlseucl} such
  that there exists $u_-\in L^2(\R^d)$ and 
  \begin{equation*}
    \left\| u(t)-\gd(t)u_-\right\|_{L^2}= \left\|
    \gd(-t)u(t)-u_-\right\|_{L^2} \to 0 \quad \text{as }t\to -\infty.
  \end{equation*}
Then $u\equiv 0$ and $u_-\equiv 0$.
\end{prop}
\begin{proof}[Sketch of the proof]
  Let $\psi \in C_0^\infty (\R^d)$ and $t_1\le t_2\le -T$: by assumption, 
  \begin{equation*}
    \< \psi, \gd(-t_2)u(t_2)-\gd(-t_1)u(t_1)\>= -i\int_{t_1}^{t_2}\<
    \gd(\tau)\psi, \(|u|^{2\si}u\)(\tau)\>d\tau
  \end{equation*}
goes to zero as $t_1,t_2\to -\infty$. But for $\tau\to -\infty$, we have
\begin{equation*}
  \gd(\tau)\psi \sim c \frac{e^{i|x|^2/(4\tau)}}{|\tau|^{d/2}}\widehat
  \psi\(\frac{x}{2\tau}\)\quad ;\quad u(\tau)\sim \gd(\tau)u_-\sim c
  \frac{e^{i|x|^2/(4\tau)}}{|\tau|^{d/2}}\widehat 
  {u_-}\(\frac{x}{2\tau}\).
\end{equation*}
Therefore,
\begin{equation*}
  \<
    \gd(\tau)\psi, \(|u|^{2\si}u\)(\tau)\> \sim
    \frac{C}{|t|^{\si d}}\< \widehat
  \psi, \left|\widehat  {u_-} \right|^{2\si}\widehat  {u_-}\>.
\end{equation*}
This function of $\tau$ is not integrable, unless
\begin{equation*}
  \< \widehat
  \psi, \left|\widehat  {u_-} \right|^{2\si}\widehat  {u_-}\>=0.
\end{equation*}
Since $\psi \in C_0^\infty (\R^d)$ is arbitrary, this means that
$\widehat  {u_-} \equiv 0\equiv u_-$. The assumption and the
conservation of mass then imply $u\equiv 0$.
\end{proof}
When $\si\le 1/d$, long range effects must be taken into account, \emph{even
in a radial setting} (see e.g. \cite{Ozawa,CaCMP} for the case $d=1$,
\cite{GO93} for $d\ge 2$). 
\subsection{Asymptotic completeness in $H^1$}
\label{sec:asymeucl}

To get a complete picture of large time behavior of solutions to
\eqref{eq:nlseucl}, we proceed to the next step which consists in
establishing asymptotic 
completeness, that is, proving that the wave operators $W_\pm$ are
invertible on their range. Here we only recall some results in $H^1(\R^d)$,
and we do not mention  what can be done in weaker Sobolev spaces, or
in weighted Sobolev spaces (see e.g. \cite{CW92,GV79Scatt,NO,VZ}).  
%We shall consider the invertibility of $W_+$,
%since the invertibility of $W_-$ is proved in a similar fashion. 
\smallbreak

The original proof of the asymptotic completeness for \eqref{eq:nlseucl}
in $H^1(\R^d)$ is due to 
Ginibre and Velo \cite{GV85}. Let $t_0=0$ and $\varphi \in H^1(\R^d)$:
the local in time $H^1$ solution to the Cauchy 
problem \eqref{eq:nlseucl}--\eqref{eq:cieucl} is actually global in
time for $0<\si<2/(d-2)$, thanks to the conservations of mass and
energy, since the nonlinearity is defocusing: 
\begin{align*}
  \|u(t)\|_{L^2}&=\|\varphi\|_{L^2},\\
\|\nabla u(t)\|_{L^2}^2 +\frac{1}{\si +1}\|u(t)\|_{L^{2\si +2}}^{2\si
  +2}& =\|\nabla \varphi\|_{L^2}^2 +\frac{1}{\si
  +1}\|\varphi\|_{L^{2\si +2}}^{2\si 
  +2}.
\end{align*}
Using Morawetz inequality and dispersive
estimates for $\gd(t)$, they prove that
\begin{equation*}
  \|u(t)\|_{L^q(\R^d)}\to 0 \quad\text{as }t\to \pm\infty,\quad \forall
  q\in \left]2,\frac{2}{d-2}\right[.
\end{equation*}
This makes it possible to show that $u\in L^p(\R;L^q(\R^d))$ for all
$d$-admissible pairs $(p,q)$, as soon as $2/d<\si<2/(d-2)$. Asymptotic
completeness follows easily:
\begin{prop}[\cite{GV85}, see also
  \cite{CazCourant}]\label{prop:CAGV} 
 Let $d\ge 3$, $t_0=0$ and $\varphi \in H^1(\R^d)$. If
 $2/d<\si<2/(d-2)$, then there exist $u_\pm\in H^1(\R^d)$ 
such that
\begin{equation*}
  \left\| u(t)-\gd(t)u_\pm\right\|_{H^1(\R^d)}\to 0\quad \text{as
  }t\to \pm\infty .
\end{equation*}
\end{prop}
More recently, a simplified
proof was proposed by Tao, Visan and Zhang \cite{TaoVisanZhang},
relying on an \emph{interaction} Morawetz inequality as introduced in
\cite{CKSTTCPAM}. We recall this approach  for essentially two reasons:
\begin{itemize}
\item It is shorter than the original one
\cite{GV85} (or \cite{CazCourant}).
\item It does not use dispersive estimates for $\gd(t)$.
\end{itemize}
The second point seems to be crucial to prove Theorem~\ref{theo:CA} as a
consequence of the proof in \cite{TaoVisanZhang} and of the
interaction Morawetz inequality that we establish in
Section~\ref{sec:morawetz}. 
The interaction Morawetz inequality presented in \cite{TaoVisanZhang}
reads as follows:
\begin{prop}[\cite{TaoVisanZhang}]\label{prop:CATVZ}
  Let $d\ge 3$, $t_0=0$ and $\varphi\in H^1(\R^d)$. Let $I$ be a
  compact time interval. There exists $C$ independent of $I$
  such that the following holds.
  \begin{itemize}
  \item If $d=3$, then the solution to
  \eqref{eq:nlseucl}--\eqref{eq:cieucl} satisfies:
    \begin{equation}\label{eq:mor3}
      \int_I\int_{\R^3}|u(t,x)|^4dxdt \le C \|u\|_{L^\infty (I;H^1)}^4.
    \end{equation}
  \item If $d\ge 4$, then the solution to
  \eqref{eq:nlseucl}--\eqref{eq:cieucl} satisfies:
    \begin{equation}\label{eq:mor4+}
      \int_I\int_{\R^3}\frac{|u(t,y)|^2|u(t,x)|^2}{|x-y|^3}dxdydt \le
      C \|u\|_{L^\infty (I;H^1)}^4.
    \end{equation}
  \end{itemize}
These inequalities imply that there exists $\widetilde C$ independent of $I$
  such that:
  \begin{equation}
    \label{eq:morgen}
    \big\| u\big\|_{L^{d+1}\(I;L^\frac{2(d+1)}{d-1}(\R^d)\)}\le \widetilde C
    \|u\|_{L^\infty (I;H^1)}. 
  \end{equation}
\end{prop}
In Section~\ref{sec:morawetz}, we establish the analogue of
\eqref{eq:mor3}--\eqref{eq:mor4+} on the hyperbolic space $\H^n$, so
we do not recall how \eqref{eq:mor3} and \eqref{eq:mor4+} are proven
here: the method on $\H^n$ is similar, with an additional drop of
geometry. 
\smallbreak

If $d=3$, then \eqref{eq:morgen} is exactly \eqref{eq:mor3}. On the
other hand, when $d\ge 4$, \eqref{eq:morgen} follows from
\eqref{eq:mor4+} by interpreting the convolution with $\frac{1}{|x|^3}$ as differentiation, and thanks to the inequality (\cite[Lemma~5.6]{TaoVisanZhang})
\begin{equation}\label{eq:paradiff}
  \left\| \left|\nabla\right|^{-\frac{d-3}{4}}f\right\|_{L^4}^2\lesssim
  \left\| \left|\nabla\right|^{-\frac{d-3}{2}}|f|^2\right\|_{L^2} ,
\end{equation}
which can be established by using paradifferential calculus. 
\smallbreak
\begin{rema}
The pair $\( d+1,\frac{2(d+1)}{d-1}\)$ present in
\eqref{eq:morgen} is $2$-admissible.  
\end{rema}
Using the \emph{a priori}
estimate provided by the conservations of mass and energy, one infers the  \emph{a
  priori} bound
\begin{equation*}
  \big\| u\big\|_{L^{d+1}\(\R;L^\frac{2(d+1)}{d-1}(\R^d)\)}\le \widetilde C
    \|u\|_{L^\infty (\R;H^1)}\lesssim \|\varphi\|_{H^1}.
\end{equation*}
Let $\eta>0$ be a small constant to be fixed later. The line $\R$ can
be divided into $J$ (for some \emph{finite} $J$ from the above estimate)
subintervals $I_j=[\tau_j,\tau_{j+1}]$ such that 
\begin{equation*}
  \big\| u\big\|_{L^{d+1}\(I_j;L^\frac{2(d+1)}{d-1}(\R^d)\)}\le \eta.
\end{equation*}
 As a consequence of \cite[Lemma~2.7]{TaoVisanZhang}, if $d\ge 3$ and
 $2/d<\si<2/(d-2)$, there exist $C>0$, $\delta \in ]0,1[$ and a
 $d$-admissible pair $(p_0,q_0)$ such that for any time interval $I$:
 \begin{equation}
   \label{eq:12h38}
   \left\| |u|^{2\si}u\right\|_{L^2\(I;W^{1,\frac{2d}{d+2}}\)} \le C 
\left\| u\right\|_{L^{d+1}\(I;L^\frac{2(d+1)}{d-1}\)}^{2\si \delta}
\left\| u\right\|_{L^\infty\(I;H^1\)}^{2\si (1-\delta)}
\left\| u\right\|_{L^{p_0}\(I;W^{1,q_0}\)}.
 \end{equation}
This estimate follows from H\"older inequality (see
\cite{TaoVisanZhang}), and algebraic computations on $d$-admissible
pairs. Using Strichartz estimates on \eqref{eq:duhameleucl} (with
$t_0$ replaced by $\tau_j$ and $\varphi$ replaced by $u(\tau_j)$), and
\eqref{eq:12h38}, we get, for $1\le j\le J$, and any $d$-admissible
pair $(p,q)$,
\begin{align*}
  \|u\|_{L^p(I_j;W^{1,q})} &\lesssim \|u(\tau_j)\|_{H^1} + \left\|
  |u|^{2\si}u\right\|_{L^2\(I;W^{1,\frac{2d}{d+2}}\)}\\
&\lesssim \|u\|_{L^\infty(\R;H^1)} + \left\|
  u\right\|_{L^{d+1}\(I_j;L^\frac{2(d+1)}{d-1}\)}^{2\si \delta} 
\left\| u\right\|_{L^\infty\(I_j;H^1\)}^{2\si (1-\delta)}
\left\| u\right\|_{L^{p_0}\(I_j;W^{1,q_0}\)}\\
&\lesssim \|u\|_{L^\infty(\R;H^1)} + \eta^{2\si \delta} 
\left\| u\right\|_{L^\infty\(\R;H^1\)}^{2\si (1-\delta)}
\left\| u\right\|_{L^{p_0}\(I_j;W^{1,q_0}\)}.
\end{align*}
Fix $(p,q)=(p_0,q_0)$: taking $\eta>0$ sufficiently small, we find
\begin{equation*}
 \|u\|_{L^{p_0}(I_j;W^{1,q_0})} \le C \|u\|_{L^\infty(\R;H^1)},\quad
 \forall 1\le j\le J,
\end{equation*}
hence $u,\nabla u\in L^{p_0}(\R;W^{1,q_0})$. We deduce $u,\nabla u\in
L^{p}(\R;W^{1,q})$ for all $d$-admissible pairs $(p,q)$. 
Asymptotic completeness is then straightforward: let $t_2\ge t_1\ge
0$. From inhomogeneous Strichartz estimates with $(p_2,q_2)=(2,\frac{2n}{n-2})$ and \eqref{eq:12h38}, we have:
\begin{align*}
  \left\| \gd(-t_2)u(t_2)-\gd(-t_1)u(t_1)\right\|_{H^1} & \lesssim
  \left\| |u|^{2\si}u
  \right\|_{L^2\([t_1,+\infty[;W^{1,\frac{2d}{d+2}}\)} \\
& \lesssim \left\| u
  \right\|_{L^{d+1}\([t_1,+\infty[;L^{\frac{2(d+1)}{d-1}}\)}^{2\si
  \delta}.
\end{align*}
Since the last term goes to zero as $t_1\to +\infty$, this proves
Proposition~\ref{prop:CAGV} for positive time. The proof for negative
time is similar.

\begin{rema}
In Proposition~\ref{prop:CAGV}, we assume that $\sigma<2/(d-2)$. Scattering
for the  $H^1$-critical case $\sigma=2/(d-2)$ was recently proved  for
$d=3$ in  \cite{CKSTT:global}, for $d=4$ in \cite{RV:global} and
finally for $d\ge 5$ in \cite{TV:global}. These results are not yet
available in  hyperbolic spaces since their proof, among other tools,
uses very subtle arguments in Fourier analysis, arguments that are not
at hand yet in $\mathbb H\sp n$. 
\end{rema}

\section{Weighted Strichartz inequalities and consequences}
\label{sec:strichartz}

The general idea is that weighted Strichartz estimates are available
on hyperbolic space $\H^n$, provided that we restrict our study to
\emph{radial} 
functions. The weight has exponential decay in
space. This decay gives us the ``usual'' Strichartz estimates
recalled in Proposition~\ref{prop:strieucl}, for $d$-admissible pairs,
for \emph{any} $d\ge n$. 
As usual for Strichartz estimates, we distinguish the case $n\ge 3$ from
the case $n=2$. The former is easier to present, and we start with
it. In $\H^n$, we denote
\begin{equation*}
  \w(r) := \(\frac{\sinh r}{r}\)^{\frac{n-1}{2}}\quad ;\quad
  U(t)=e^{it\Delta_{\H^n}}. 
\end{equation*}

\subsection{Case $n\ge 3$}
\label{sec:stri3+}

The following global result was established in \cite{BaHyper} for $n=3$, and in
\cite{VittoriaDR} for $n\ge 4$:
\begin{prop}[Weighted Strichartz estimates in
  $\H^n$, $n\ge 3$]\label{prop:weightedstri3+} Let $n\ge 3$.\\ 
 {\rm1.} For any $n$-admissible pair $(p,q)$, there exists $C_q$
    such that
\begin{equation*}
    \left\|\w^{1-\frac{2}{q}} U(\cdot)\phi
    \right\|_{L^p({\R};L^q)}\le C_q \|\phi 
    \|_{L^2}
  \end{equation*}
 for every \emph{radial} function $\phi \in L^2_{\rm rad}({\H}^n)$.%\\

\noindent{\rm2.} For any $n$-admissible pairs $(p_1,q_1)$ and $
    (p_2,q_2)$ and any 
    interval $I$, there exists $C_{q_1,q_2}$ independent of $I$ such that 
\begin{equation*}
      \left\|\w^{1-\frac{2}{q_1}} \int_{I\cap\{s\le
      t\}} U(t-s)F(s)ds 
      \right\|_{L^{p_1}(I;L^{q_1})}\le C_{q_1,q_2} \left\|
      \w^{1-\frac{2}{q'_2}}F\right\|_{L^{p'_2}\(I;L^{q'_2}\)}
    \end{equation*}
for every \emph{radial} function $F\in L^{p'_2}\(I;L^{q'_2}_{\rm
  rad}(\H^n)\)$.  
\end{prop}
\begin{cor}\label{cor:stri3+}
  Let $d\ge n\ge 3$. Then Strichartz estimates hold for $d$-admissible
  pairs and radial functions on $\H^n$:\\
{\rm1.} For any $d$-admissible pair $(p,q)$, there exists $C_q=C_q(n,d)$
    such that
\begin{equation*}
    \left\| U(\cdot)\phi \right\|_{L^p({\R};L^q)}\le C_q \|\phi
    \|_{L^2},\quad \forall \phi \in L^2_{\rm rad}({\H}^n).
  \end{equation*}

\noindent{\rm2.} For any $d$-admissible pairs $(p_1,q_1)$ and $
    (p_2,q_2)$ and any 
    interval $I$, there exists $C_{q_1,q_2}=C_{q_1,q_2}(n,d)$
    independent of $I$ such that  
\begin{equation*}
      \left\| \int_{I\cap\{s\le
      t\}} U(t-s)F(s)ds 
      \right\|_{L^{p_1}(I;L^{q_1})}\le C_{q_1,q_2} \left\|
      F\right\|_{L^{p'_2}\(I;L^{q'_2}\)},
    \end{equation*}
for every $F\in L^{p'_2}\(I;L^{q'_2}_{\rm rad}({\H}^n)\)$.
\end{cor}
\begin{proof}
  To prove the first estimate, it is enough to prove it for the endpoint
  estimate, $(p,q)=(2,\frac{2d}{d-2})$. Define $s$ by
  \begin{equation*}
    \frac{1}{n}=\frac{1}{d}+\frac{1}{s}\cdot
  \end{equation*}
We have $s\ge 0$, since $d\ge n$. Let $\phi\in L^2_{\rm
  rad}(\H^n)$. H\"older inequality and the first 
part of Proposition~\ref{prop:weightedstri3+} yield:
\begin{equation}\label{passagend}
  \left\| U(\cdot )\phi\right\|_{L^2\(\R;L^{\frac{2d}{d-2}}\)} \le \left\|
  \w^{2/n}U(\cdot )\phi\right\|_{L^2\(\R;L^{\frac{2n}{n-2}}\)} \left\|
  \w^{-2/n}\right\|_{L^s}
\end{equation}
  $$\le C_{\frac{2n}{n-2}}
  \left\|\phi\right\|_{L^2(\H^n)} \left\| 
  \w^{-2/n}\right\|_{L^s} . $$

If $d=n$, then $s=\infty$, and we have obviously $\w^{-2/n} \in
L^\infty$. If $d>n$, then $\w^{-2/n}\in L^s$ if and only if
\begin{equation*}
  \int_0^\infty \( \frac{r}{\sinh
  r}\)^{s\frac{n-1}{n}}\(\sinh r\)^{n-1}dr <\infty.
\end{equation*}
This integral is convergent, since $s>n$ ($d$ is finite). The first
estimate of the corollary follows by interpolation, by conservation of
the $L^2$ norm.
\smallbreak 

We turn to the inhomogeneous estimates. Let $(p_1,q_1)$ and
$(p_2,q_2)$ be $d$-admissible pairs. Let $(p_1,r_1)$ and
$(p_2,r_2)$ be the corresponding $n$-admissible pairs:
\begin{equation*}
  \frac{2}{p_j}=d\( \frac{1}{2}-\frac{1}{q_j}\)=n\(
  \frac{1}{2}-\frac{1}{r_j}\).
\end{equation*}
Note that since $d\ge n$, $q_j \le r_j$. Therefore,
$s_j$, given by 
\begin{equation*}
  \frac{1}{q_j}=\frac{1}{r_j}+\frac{1}{s_j},\virgp
\end{equation*}
is non-negative. Using H\"older inequality and the second part of
Proposition~\ref{prop:weightedstri3+}, we find:
\begin{align*}
  \left\|\int_{I\cap\{s\le
      t\}} U(t-s)F(s)ds 
      \right\|_{L^{p_1}(I;L^{q_1})} & \le\\
      \le \Big\|\w^{1-\frac{2}{r_1}}\int_{I\cap\{s\le 
      t\}}& U(t-s)F(s)ds 
      \Big\|_{L^{p_1}(I;L^{r_1})}
      \left\|\w^{-1+\frac{2}{r_1}}\right\|_{L^{s_1}}\\
& \lesssim 
      \left\|\w^{1-\frac{2}{r_2'}}F 
      \right\|_{L^{p_2'}(I;L^{r_2'})}
      \left\|\w^{-1+\frac{2}{r_1}}\right\|_{L^{s_1}}\\
& \lesssim 
      \left\|F 
      \right\|_{L^{p_2'}(I;L^{q_2'})}
\left\|\w^{1-\frac{2}{r_2'}}\right\|_{L^{s_2}}
      \left\|\w^{-1+\frac{2}{r_1}}\right\|_{L^{s_1}}.
\end{align*}
Therefore, we have to check that 
$\w^{-1+\frac{2}{r_j}}\in   L^{s_j}(\H^n)$.
If $q_j=2$, then $r_j=2$ and $s_j=\infty$. If $q_j>2$, then
the above integrability condition is equivalent to:
\begin{equation*}
  s_j\(\frac{1}{2}-\frac{1}{r_j}\) > 1
  \Leftrightarrow \frac{1}{2}-\frac{1}{r_j}>
  \frac{1}{s_j} = \frac{1}{q_j}-\frac{1}{r_j}.
\end{equation*}
Since $q_j>2$, this is satisfied, and the corollary follows.
\end{proof}

\subsection{Case $n=2$}
\label{sec:stri2}
When $n=2$, the analogue of Proposition~\ref{prop:weightedstri3+} is
not proven, but we have from \cite{BaHyper}:
\begin{equation}\label{eq:solfond2d}
  e^{it\Delta_{\H^2}}\phi (\Omega)=
  \frac{c}{|t|^{3/2}}e^{-it/2}\int_{\H^2} \phi(\Omega')\int_\rho^\infty
  \frac{se^{is^2/4t}}{\sqrt{\cosh s -\cosh\rho }}dsd\Omega', 
\end{equation}
where $\rho = d(\Omega,\Omega')$. 
The following weighted dispersion estimate holds for \emph{radial} functions
  in $\H^2$. Denote 
  \begin{equation*}
    \widetilde {\tt w}_2(r)=\(\frac{\sinh r}{r(1+r)}\)^{1/2}\cdot
  \end{equation*}
\begin{prop}\label{prop:disp2d}
  Let $\eps\in ]0,1[$. There exists $C_\eps>0$ such
  that  
  \begin{equation*}
  \widetilde {\tt w}_2^{1-\eps}(r)  e^{it\Delta_{\H^2}}\phi(\Omega)
    \le  \frac{C_\eps}{|t|^{3/2}}\int_{\H^2}|\phi(\Omega')|
    \frac{d\Omega'}{\widetilde {\tt w}_2^{1-\eps}(r')},\quad \forall t\not = 0,
    \quad \forall\phi\in L^1_{\rm rad}(\H^2),
  \end{equation*}
where $r=d(0,\Omega)$ and $r'=d(0,\Omega')$.
\end{prop}
\begin{rema}
  For small time, this weighted estimate is worse than the one in $\R^2$ in
  terms of powers of $t$: $|t|^{-3/2}$ instead of
  $|t|^{-1}$. Formally, we could get a rate in $|t|^{-1-\eps}$ by
  integration by parts in the $s$ integral in \eqref{eq:solfond2d},
  by considering derivatives of order $1/2 -\eps$. We shall not
  pursue this approach here, and content ourselves with
  Proposition~\ref{prop:disp2d}. 
\end{rema}
\begin{proof}
We first prove
  \begin{equation}\label{eq:decayHeat}
    \left|\int_\rho^\infty
  \frac{se^{is^2/4t}}{\sqrt{\cosh s -\cosh\rho }}ds\right|\le
  C\(\frac{\rho}{\sinh \rho}\)^{1/2}\sqrt{1+\rho},
  \end{equation}
where $C$ is independent of $\rho\ge 0$. Note that this estimate is
analogous to the one given in \cite{DaviesHeat}: for the heat
operator, an additional Gaussian decay is available (replace
$e^{is^2/4t}$ with $e^{-s^2/4t}$). The computation
below shows that this extra decay   is not necessary in order for \eqref{eq:decayHeat} to be
true. We have obviously
\begin{equation*}
 \left|\int_\rho^\infty
  \frac{se^{is^2/4t}}{\sqrt{\cosh s -\cosh\rho }}ds\right| \le \int_\rho^\infty
  \frac{s}{\sqrt{\cosh s -\cosh\rho }}ds .
\end{equation*}
Using "trigonometry", we find:
\begin{equation*}
  \int_\rho^\infty
  \frac{s}{\sqrt{\cosh s -\cosh\rho }}ds = \int_\rho^\infty
  \frac{s}{\sqrt{2\sinh\(\frac{s+\rho}{2}\)\sinh\(\frac{s-\rho}{2}\)
  }}ds.  
\end{equation*}
With the change of variable $y=s-\rho$, we estimate:
\begin{equation*}
  \int_0^\infty \frac{y+\rho}{\sqrt{\sinh\(\rho+\frac{y}{2}\)
\sinh\(\frac{y}{2}\)
  }}dy=\int_0^\infty \frac{(y+2\rho)-\rho}{\sqrt{\sinh\(\rho+\frac{y}{2}\)
\sinh\(\frac{y}{2}\)
  }}dy .
\end{equation*}
For the first term, we use the fact that
\begin{equation*}
  s\mapsto \frac{s}{\sinh s} \text{ is non-increasing},
\end{equation*}
to have the estimate:
\begin{align*}
  \int_0^\infty \frac{y+2\rho}{\sqrt{\sinh\(\rho+\frac{y}{2}\)
\sinh\(\frac{y}{2}\)
  }}dy &\le \sqrt{2\frac{\rho}{\sinh \rho}}\int_0^\infty
\frac{\sqrt{y+2\rho}}{\sqrt{\sinh\(\frac{y}{2}\)
  }}dy\\
\lesssim \sqrt{\frac{\rho}{\sinh \rho}}&\(\int_0^\infty
\sqrt{\frac{y}{\sinh\(\frac{y}{2}\) }}dy + \sqrt\rho \int_0^\infty
\frac{dy}{\sqrt{\sinh\(\frac{y}{2}\) }} \)\\
\lesssim \sqrt{\frac{\rho}{\sinh \rho}}&\sqrt{1+\rho}.
\end{align*}
For the second term, we have:
\begin{align*}
  \rho \int_0^\infty \frac{1}{\sqrt{\sinh\(\rho+\frac{y}{2}\)
\sinh\(\frac{y}{2}\)  }}dy \le \frac{\rho}{\sqrt{\sinh \rho}}\int_0^\infty
\frac{dy}{\sqrt{\sinh\(\frac{y}{2}\) }},
\end{align*}
and \eqref{eq:decayHeat} follows. To infer the proposition, we mimic the
computations of \cite{BaHyper}, \S5. From \eqref{eq:solfond2d} and
\eqref{eq:decayHeat}, we have:
\begin{equation*}
   \left|e^{it\Delta_{\H^2}}\phi (\Omega)\right| \lesssim
  \frac{1}{|t|^{3/2}}\int_{\H^2}\left| \phi(\Omega')\right|
  \(\frac{\rho}{\sinh \rho}\)^{1/2}\sqrt{1+\rho} \,d\Omega'. 
\end{equation*}
Recall that $\phi$ is radial: with the usual abuse of notations,
\begin{equation*}
  \phi(\Omega') = \phi(\cosh r',\om'\sinh r')= \phi(r').
\end{equation*}
Using hyperbolic coordinates, 
\begin{equation*}
  \rho = d(\Omega,\Omega')= \cosh^{-1}\(\cosh r\cosh r' - \sinh
  r\sinh r' \om\cdot \om'\),
\end{equation*}
and we can write:
\begin{equation*}
  \left|e^{it\Delta_{\H^2}}\phi (\Omega)\right| \lesssim
  \frac{1}{|t|^{3/2}}\int_{\H^2}\left| \phi(\Omega')\right|
  K(r,r') d\Omega',
\end{equation*}
where the kernel $K$ is given by:
\begin{equation*}
  K(r,r')= \int_{\Sph^1}\sqrt{f\(\cosh^{-1}\(\cosh r \cosh r' - \sinh
  r\sinh r' \om\cdot \om'\)\)}d\om',
\end{equation*}
with
\begin{equation*}
  f(y)=\frac{y}{\sinh y}(1+y). 
\end{equation*}
With $x= \om\cdot \om'$, we have:
\begin{equation*}
  K(r,r')= \int_{-1}^1\sqrt{f\(\cosh^{-1}\(\cosh r \cosh r' - x\sinh
  r\sinh r' \)\)}\frac{dx}{\sqrt{1-x^2}}.
\end{equation*}
The lemma follows from H\"older inequality. For $\l>0$, we have
\begin{equation*}
  K(r,r')\le C_\l\( \int_{-1}^1f\(\cosh^{-1}\(\cosh r \cosh r' - x\sinh
  r\sinh r' \)\)^{1+\l}dx\)^{1/(2+2\l)}.
\end{equation*}
With the change of variable
\begin{equation*}
  \cosh y = \cosh r \cosh r' - x\sinh
  r\sinh r', 
\end{equation*}
this yields
\begin{align*}
  K(r,r')&\le C_\l \( \frac{1}{\sinh r\sinh r'}
\int_{|r-r'|}^{r+r'}f(y)^{1+\l}\sinh y
dy\)^{1/(2+2\l)}\\
&\le C'_\l\( \frac{1}{\sinh r\sinh r'}
\int_{|r-r'|}^{r+r'}y(1+y)dy\)^{1/(2+2\l)}\\
&\le C_\l' \( \frac{(r+r')^3-|r-r'|^3}{\sinh r\sinh
  r'}\)^{1/(2+2\l)} \\
&\le C_\l'' \( \frac{r'(1+r')r(1+r)}{\sinh r\sinh
  r'}\)^{1/(2+2\l)} . 
\end{align*}
This completes the proof of the proposition, with
$\frac{1}{1+\l}=1-\eps$. 
\end{proof}
We find a weighted dispersive estimate and weighted Strichartz
estimates which are
similar to the ones for
radial functions in $\H^3$. The difference is that we must replace
${\tt w}_3$ with $\widetilde {\tt w}_2^{1-\eps}$. Even though the value
$\eps=0$ is excluded, we can consider $d>3$ arbitrarily close to
$3$ and repeat the argument in (\ref{passagend}). On the other hand, since $\widetilde {\tt w}_2$ is bounded, we
have the  Strichartz estimates as in $\R^3$ for free: this yields the
Corollary~\ref{cor:stri3+} with $n=2$ and $d\ge 3$.\par
>From Proposition \ref{prop:disp2d} we have the dispersion
\begin{equation*}
 \left\|e^{it\Delta_{\H^2}}\phi
 \right\|_{L^\infty_{\rm rad}(\H^2)}\lesssim
 \frac{1}{|t|^{3/2}}\|\phi\|_{L^1_{\rm rad}(\H^2)} ,\quad \forall
 t\not =0.
\end{equation*}
For small time, \cite[Theorem~1.2]{BaHyper} yields:
\begin{equation*}
 \left\|e^{it\Delta_{\H^2}}\phi
 \right\|_{L^\infty(\H^2)}\lesssim
 \frac{1}{|t|}\|\phi\|_{L^1(\H^2)} ,\quad \forall
 t\in [-1,1]\setminus\{0\}.
\end{equation*}
We infer the global dispersive estimate, for $2\le d\le 3$ (not
necessarily an integer):
\begin{equation*}
 \left\|e^{it\Delta_{\H^2}}\phi
 \right\|_{L^\infty_{\rm rad}(\H^2)}\le
 \frac{C_d}{|t|^{d/2}}\|\phi\|_{L^1_{\rm rad}(\H^2)} ,\quad \forall
 t\not =0.
\end{equation*}
We conclude:
\begin{cor}\label{cor:stri2}
  Let $d\ge 2$. Then Strichartz estimates hold for $d$-admissible
  pairs and radial functions on $\H^2$:\\
{\rm1.} For any $d$-admissible pair $(p,q)$, there exists $C_q=C_q(d)$
    such that
\begin{equation*}
    \left\| U(\cdot)\phi \right\|_{L^p({\R};L^q)}\le C_q \|\phi
    \|_{L^2},\quad \forall \phi \in L^2_{\rm rad}(\H^2).
\end{equation*}

\noindent{\rm2.} For any $d$-admissible pairs $(p_1,q_1)$ and $
    (p_2,q_2)$ and any 
    interval $I$, there exists $C_{q_1,q_2}=C_{q_1,q_2}(d)$
    independent of $I$ such that  
\begin{equation*}
      \left\| \int_{I\cap\{s\le
      t\}} U(t-s)F(s)ds 
      \right\|_{L^{p_1}(I;L^{q_1})}\le C_{q_1,q_2} \left\|
      F\right\|_{L^{p'_2}\(I;L^{q'_2}\)},
\end{equation*}
for every $F\in L^{p'_2}\(I;L^{q'_2}_{\rm rad}({\H}^2)\)$.
\end{cor}

\section{Scattering for small data in $L^2_{\rm rad}(\H^n)$}
\label{sec:small}

The proof of Theorem~\ref{theo:smallL2} is straightforward in view of
Section~\ref{sec:smalleucl} and Corollaries~\ref{cor:stri3+} and
\ref{cor:stri2}. 
\smallbreak

Let $n\ge 2$ and $0<\si\le 2/n$. Set $d = 2/\si$:
Corollaries~\ref{cor:stri3+} and 
\ref{cor:stri2} yield the same Strichartz estimates as in
$\R^d$, provided that we work with radial functions. Simply notice
that the proof of Proposition~\ref{prop:L2eucl} relies only on
functional analysis: H\"older inequality and Strichartz
estimates. Theorem~\ref{theo:smallL2} follows: the statement is the
analogue of Proposition~\ref{prop:L2eucl}, with $d = 2/\si$.
\begin{rema}
 This result shows that  there are
   no long range effects, at least in a neighborhood of the origin in
   $L^2_{\rm rad}(\H^n)$. This can be compared to
   \cite[Proposition~1.1]{CaDCDS}. There, the following nonlinear
   Schr\"odinger equation is considered:
   \begin{equation*}
     i\d_t u +\frac{1}{2}\Delta_{\R^n}u = -\frac{x_1^2}{2}u +
     V(x_2,\ldots, x_n)u +\kappa |u|^{2\si}u \quad ; \quad x\in \R^n,\
     n\ge 1, 
   \end{equation*}
where $V$ is any quadratic polynomial ($V\equiv 0$ if $n=1$). It is
proved that for $0<\si\le 2/n$, there is a small data scattering
theory in $L^2$, just as in Theorem~\ref{theo:smallL2}. This is 
because the repulsive potential $-x_1^2$ yields an 
exponential decay in time of the free solution. Here, this exponential
decay in time is replaced by an exponential decay in space. The proof
relies on the same idea though: we have Strichartz estimates that make
it possible to pretend that we work in $\R^d$ with $d\ge n$.     
\end{rema}

\section{Wave operators in $H^1_{\rm rad}(\H^n)$}
\label{sec:wave}

The argument for the proof of Theorem~\ref{theo:waveH1} is similar.
First, taking 
$d=n$, we cover the range $2/n\le \si<2/(n-2)$. To cover the range
$0<\si<2/n$, we keep the
value $d=2/\si$: $d>n$. The proof of Proposition~\ref{prop:waveH1eucl} uses
the same arguments as the proof of Proposition~\ref{prop:L2eucl},
plus Sobolev embeddings. Therefore, we simply have to check that the
step where Sobolev embeddings are used can be adapted.
\smallbreak

Recall that we work in $\H^n$, and that we pretend that we work in
$\R^d$, with $d=2/\si\ge n$. In the proof of
Proposition~\ref{prop:waveH1eucl}, we used the embedding:
\begin{equation*}
  H^1(\R^d)\subset L^{2\si+2}(\R^d). 
\end{equation*}
Since we assume $\si<2/(n-2)$, we have:
\begin{equation*}
  H^1(\H^n)\subset L^{2\si+2}(\H^n). 
\end{equation*}
Therefore, we can argue as in Section~\ref{sec:small}: we can mimic
the approach to prove Proposition~\ref{prop:waveH1eucl}, which is
based on functional analysis.

\section{Morawetz estimates in $\H^n$}
\label{sec:morawetz}
In this section we prove some  Morawetz type estimates for general
solutions to the equation  
\eqref{eq:nls}. Note that in this section  the solutions are not
necessarily radial. We start by stating these inequalities. Their
proof will be a consequence of more general geometric set up already
used by Hassell, Tao and Wunsch \cite{HTW05}  while studying the same type of
estimates for non-trapping asymptotically conic manifolds. This last
paper in turn was based on the interaction Morawetz inequality
introduced by Colliander et al. \cite{CKSTTCPAM}.

Define the operator $ \tilde H$ acting on a function 
$f$ defined on $\H^n\times \H^n$ as 
\begin{equation*}
  \tilde H f(\Omega,\Omega')=-\Delta_{\H^n\times \H^n}f(\Omega,\Omega')=-\Delta_{\Omega}f(\Omega,\Omega')-\Delta_{\Omega'}f(\Omega,\Omega'),
\end{equation*}
and the distance function $d(\Omega,\Omega')={\rm
  dist}_{\H^n}(\Omega,\Omega')$. 
Then we can state the following theorem
\begin{theo}\label{mor:gen} 
For any compact interval of time  $I$ and for any $u$ solution to
\eqref{eq:nls}, 
\begin{equation*}%\label{mor:geneq}
-\int_I\int_{\H^n\times \H^n}\tilde
H^2(d(\Omega,\Omega'))|u(t,\Omega)|^2|u(t,\Omega')|^2\, 
\, d\Omega  d\Omega'\, dt\le
C\|u\|_{L^\infty (I;H^1)}^4. 
\end{equation*}
\end{theo}
Since
\begin{equation*}
-\tilde H^2\(d(\Omega,\Omega')\)=\left\{
  \begin{aligned}
    \delta_{\Omega'}(\Omega)+\delta_{\Omega}(\Omega') &\quad\mbox{ if }  n=3,\\
 2\frac{\cosh(d(\Omega,\Omega'))}{\sinh^3(d(\Omega,\Omega'))}&\quad\mbox{
 if }  n>3, 
  \end{aligned}
\right.
\end{equation*}
we also have the following corollary:
\begin{cor}\label{cor:morHn}
For any compact interval of time  $I$ and for any $u$ solution to
\eqref{eq:nls} we have: 
\begin{itemize}
\item If $n=3$:
  \begin{equation}
    \label{mor:n=3}
    \int_I\int_{\H^3}|u(t,\Omega)|^4\, d\Omega\, dt
\le C\|u\|_{L^\infty (I;
H^1)}^4.
  \end{equation}
\item If $n>3$:
  \begin{equation*}
    %\label{mor:n>3}
    \int_I
\int_{\H^n}\int_{\H^n}\frac{\cosh(d(\Omega,\Omega'))}{\sinh^3(d(\Omega,\Omega'))} 
  |u(t,\Omega')|^2  
|u(t,\Omega)|^2\,d\Omega\,d\Omega'\,dt
\le
C\|u\|_{L^\infty (I;H^1)}^4.
  \end{equation*}
\end{itemize}
\end{cor}
\begin{rema}\label{rema1}
 Few comments are in order at this  point. First we note that in the
 above results, we  
 do not assume that $u$ is radial.  Second we recall that while the
 Morawetz type estimate proved in  \cite{HTW05} was local in time,
 ours is global, just like in the Euclidean case in \cite{CKSTTCPAM}.
 More comments will be made about this fact at the end of the proof of
 Theorem \ref{mor:gen} and Corollary \ref{cor:morHn}.

\end{rema}
Let $M$ be a general Riemannian manifold with metric $g$. We denote by
$\<\cdot,\cdot\>_g$ the product on the tangent space given by the metric
$g$. We define the real inner product for functions on $M$ 
\begin{equation}\label{prodM}
\langle u,v\rangle_M=\RE\int_M u(z)\overline{v}(z)\,dV_g(z).
\end{equation}
We will often use the commutator $[A,B]$ among pseudo-differential
operators $A$ and $B$ defined as  $[A,B]=AB-BA$. We have the following
lemma corresponding to \cite[Lemma~2.1]{HTW05}. 
\begin{lem}
Let $a(x)$ be  a real-valued tempered distribution on  a manifold $M$,
acting as  a multiplier operator $(af)(x)=a(x)f(x)$ on Schwartz
functions. Then we have the commutator identities, with $H=-\Delta_M$:
\begin{equation}\label{identity1}
i[H,a]=-i\langle \nabla a,\nabla\rangle_g+iHa=-i(\nabla^\alpha
a)\nabla_\alpha +iH_a ,
\end{equation}
and the double commutator identity
\begin{equation}\label{identity2}
-[H,[H,a]]=-\nabla_\beta {\rm Hess}(a)^{\alpha\beta}\nabla_\alpha-(H^2a),
\end{equation}
where ${\rm Hess}(a)^{\alpha\beta}$ is the symmetric tensor
\[{\rm Hess}(a)^{\alpha\beta}=(\nabla d
a)^{\alpha\beta}=g^{\alpha\gamma}g^{\beta\delta}(
\partial_\gamma\partial_\delta+\Gamma^\rho_{\gamma\delta}\partial_\rho 
a).\] 
\end{lem}
We now assume that $U$ is solution to 
\begin{equation}\label{U}
  i\d_t U +\Delta_M U =F.
\end{equation}
Then it is easy
to see that given a real pseudo-differential operator $A$ on $M$ we
have  
\begin{eqnarray}
\partial_t\langle A\, U(t),U(t)\rangle_M&=&\langle
i[H,A]U(t),U(t)\rangle_M+\label{derivative}\\
& &+\langle -iA\, F(t), U(t)\rangle_M+\langle iA\, U(t),
F(t)\rangle_M\nonumber.  
\end{eqnarray}
Next, given a real valued tempered distribution $a$ and a function
$U(t,x)$ we define 
\begin{equation}\label{Ma}
M_a(t)=\langle i[H,a]U(t),U(t)\rangle_M.
\end{equation}
By using \eqref{identity1} and the definition \eqref{prodM} we recover the familiar form of the first order momentum 
\begin{equation}\label{momentum}
M_a(t)=\IM\int_M\langle \nabla a,\nabla U(t)\rangle_g\overline{U}(t)dV_g.
\end{equation}
Now, by taking $U$ solution to 
\eqref{U} and  $A=i[H,a]$ in \eqref{derivative}  and using
\eqref{identity2}, one gets 
\begin{eqnarray}
\frac{d}{dt}\,M_a(t)&=&-\langle [H,[H,a]]U(t),U(t)\rangle_M+\label{dtMa}\\
& &+\langle [H,a]\, F(t), U(t)\rangle_M+\langle -[H,a]\, U(t),
F(t)\rangle_M\nonumber\\ 
&=&-\langle H^2a\, U(t),U(t)\rangle_M-\langle \nabla_\beta
{\rm Hess}(a)^{\alpha\beta}\nabla_\alpha U(t),U(t)\rangle_M
\nonumber\\ 
& &+\langle [H,a]\, F(t), U(t)\rangle_M+\langle -[H,a]\, U(t),
F(t)\rangle_M. \nonumber 
\end{eqnarray}

\begin{lem}
If $M$ is a Riemannian manifold with a non-positive sectional curvature
and if $a$ is a distance function defined on $M$, that is $|\nabla
a|=1$, then for any smooth function $\phi$,
\begin{equation}\label{hessian}
\langle {\rm Hess}(a)^{\alpha\beta}
\nabla_\alpha \phi,\nabla_\beta \phi\rangle_g\ge 0.
\end{equation}
\end{lem}
\begin{proof}
This is a well-known result in Riemannian geometry. We refer the reader
for example to Theorem~3.6 in \cite{Petersen}.
\end{proof}
Using \eqref{dtMa} and \eqref{hessian}  after an integration by parts
in space variable, 
we obtain for all $T\ge 0$  the key inequality: 
\begin{eqnarray}
M_a(T)-M_a(0)\ge \int_0^T-\langle H^2a\, U(t),U(t)\rangle_M+\label{key}\\
+\int_0^T\langle [H,a]\, F(t), U(t)\rangle_M+\langle -[H,a]\, U(t),
F(t)\rangle_M\, dt \nonumber. 
\end{eqnarray}
We are now ready to prove Theorem~\ref{mor:gen}.
\begin{proof}[Proof of Theorem~\ref{mor:gen}]
Following again the argument in \cite{HTW05}, we assume now that
$M=\H^n\times \H^n$, with the usual metric $\tilde g=g\otimes g$. Assume
also that $u$ is a solution to the equation \eqref{eq:nls}. It is easy
to show that   $U(t,\Omega,\Omega'):=u(t,\Omega)\,u(t,\Omega')$ is
solution to the equation 
\begin{align*}
i\partial_t U(t,\Omega,\Omega')+\Delta_\Omega\oplus \Delta_{\Omega'}
U(t,\Omega,\Omega')=& 
\(|u|^{2\si} u\)(t,\Omega)\, u(t,\Omega')+\\
&+\(|u|^{2\si} u\)(t,\Omega')\,
u(t,\Omega)\\
=:&F(t,\Omega,\Omega'). 
\end{align*}
We now set $\tilde H=-\(\Delta_\Omega\oplus \Delta_{\Omega'}\) $ and
$a(\Omega,\Omega')= 
{\rm dist}_{\H^n}(\Omega,\Omega')$.  It is easy to see that this function $a$
is a distance function with respect to the manifold $(M,
\tilde g)$. Also one can check that this manifold has  a nonpositive
sectional curvature.  Finally one can also prove, using
\eqref{identity1} and the definition of the real inner product
\eqref{prodM}, that  
\begin{eqnarray}
\langle [\tilde H,a]\, F(t), U(t)\rangle_M+\langle -[\tilde H,a]\,
U(t), F(t)\rangle_M=\\ 
=\langle -\langle\nabla a,\nabla F(t)\rangle_{\tilde g},
U(t)\rangle_M+\langle \langle\nabla a,\nabla U(t)\rangle_{\tilde g},
F(t)\rangle_M=\\ 
=2\frac{\si}{\si+1}\int_M\Delta_\Omega
a(\Omega,\Omega')\,|u|^{2\si+2}(t,\Omega)\,|u|^2(t,\Omega'). \nonumber 
\end{eqnarray}
Since (\S 5.7 of \cite{DaviesHeat})
$$\Delta_\Omega a(\Omega,\Omega')=(n-1)\coth d(\Omega,\Omega'),$$
we infer by \eqref{key} that 
\begin{equation}\label{key2}
M_a(T)-M_a(0)\ge \int_0^T-\langle \tilde H^2aU(t),U(t)\rangle_M\, dt.
\end{equation}
Using \eqref{momentum} and the fact that $a$ is a distance
function, we also have,   for all $t\in I= [0,T]$, 
\begin{equation}\label{key3}
|M_a(t)|\le C\|u\|_{L^\infty (I;H^1)}^4.
\end{equation}
The proof of the theorem now follows by combining \eqref{key2} and
\eqref{key3}.   
\end{proof}
\begin{rema}
To follow up on Remark \ref{rema1}, we can now add that while in the
above proof we were  inspired by  
\cite{HTW05},  we  obtained a global estimate thanks to the fact that
we could pick the function $a$  to be everywhere the distance between
two points, just like in the Euclidean  space. This was not possible
in \cite{HTW05}, due to the presence of ``asymptotic cones'', and as a
consequence the distance function was only good  inside a large
ball. Surprisingly enough, the distance function is not longer a good
function even in $\R^d$, for $d=1,2$. In fact in \cite{FangGrillakis},
where the case  $\R^2$ is considered,  a space localization is also
needed and again, as a consequence, the Morawetz type estimate obtained
is only local in time. 
\end{rema}

\section{Asymptotic completeness in $H^1(\H^3)$}
\label{sec:completeness}
First, we note that the proof recalled in Section~\ref{sec:asymeucl}
can be mimicked on $\H^3$: the Morawetz estimate \eqref{eq:mor3} was
adapted to the hyperbolic case, \eqref{mor:n=3}. This proves
Theorem~\ref{theo:CA} for $2/3<\si<2$. The reason why we do not have
to assume that $u$ is radial at this stage is that global in time
Strichartz estimates are available on $\H^3$, from \cite{BaHyper}. 
\smallbreak

To decrease $\si$, we proceed with the same idea as before, and
pretend that we work on $\R^d$, for $d>3$. Fix $d>3$: we first claim
that \eqref{mor:n=3} implies
\begin{equation}\label{eq:13h57}
  \big\| u\big\|_{L^{d+1}\(\R;L^\frac{2(d+1)}{d-1}(\H^3)\)}\lesssim
    \|u\|_{L^\infty (\R;H^1)}. 
\end{equation}
Indeed, we noticed in Section~\ref{sec:asymeucl} that for any $d\ge
3$, the pair $\( d+1,\frac{2(d+1)}{d-1}\)$  is
$2$-admissible. Interpolating between the pairs $(4,4)$ and
$(\infty,2)$ yields the pair $\( d+1,\frac{2(d+1)}{d-1}\)$, hence
\eqref{eq:13h57}. Having the analogue of \eqref{eq:morgen}, we can go
on with the proof presented in Section~\ref{sec:asymeucl}, thanks to
Corollary~\ref{cor:stri3+} for \emph{radial} solutions. This
completes the proof of Theorem~\ref{theo:CA} for any $0<\si<2$. Note
also that if Corollary~\ref{cor:stri3+} holds for solutions that are
not necessarily radial, then Theorem~\ref{theo:CA} will be true
without the symmetry hypothesis. 
\smallbreak

The reason why we stated Theorem~\ref{theo:CA} in the case $n=3$ only
is the following. In Section~\ref{sec:morawetz}, we have proved the
analogue of \eqref{eq:mor3} and \eqref{eq:mor4+} for the solutions to
\eqref{eq:nls} (not necessarily radial). To go on with the proof of
\cite{TaoVisanZhang}, we would need the analogue of
\eqref{eq:paradiff} on hyperbolic space $\H^n$, $n\ge 4$. A
paradifferential calculus on $\H^n$ would be welcome then, which we do
not have at hand. 
\smallbreak

If the recent proof of Tao, Visan and Zhang \cite{TaoVisanZhang}
cannot be used on $\H^n$ 
for $n\ge 4$, one might want to use the original proof of Ginibre and
Velo \cite{GV85} (or \cite{CazCourant}). As we recalled in
Section~\ref{sec:asymeucl}, this proof uses dispersive estimates for
the free Schr\"odinger group. Unfortunately, the dispersive estimates
of the free Schr\"odinger group on $\H^n$ are only \emph{local in
  time} as soon as $n\ge 4$ \cite{BaHyper}. For this
reason, even proving the same scattering results on $\H^n$ as on
$\R^n$, $n\ge 4$, does not seem obvious at all.

\appendix

\section{On the notion of criticality}
\label{sec:crit}

Consider the equation \eqref{eq:nls}. The critical scaling in $\R^n$
is 
\begin{equation*}
  s_c = \frac{n}{2} - \frac{1}{\si}\, .
\end{equation*}
At first glance, it is not clear whether the critical indices for
\eqref{eq:nls} are the same as in the Euclidean case, since the linear
part of the equation is not scale invariant. 

Note that in the case of a positive curvature, the geometry
changes the notion of criticality (see e.g. \cite{BGT,BGTMRL},
\cite{BaJMPA}), but always in the ``same order'': the positive
curvature ``creates'' more instabilities.

However, it is established in \cite{BaHyper} that there is local
well-posedness in $H^s(\H^n)$ for 
\eqref{eq:nls} if $s>s_c$: 
\begin{equation*}
  \text{subcritical in the Euclidean case}\Rightarrow
  \text{subcritical in the hyperbolic case}.
\end{equation*}
This stems from the fact that we have the same local Strichartz 
inequalities as in the Euclidean space. Moreover, in the focusing setting, 
blow-up may occur ``as in the Euclidean case'', thanks to a new virial
identity, where the negative curvature of  the hyperbolic space shows
up. 

On the other hand, some proofs of ill-posedness rely on highly
concentrated initial data and solutions, so that the geometry is not
relevant. To prove that the notion of criticality is the same in the Euclidean
and in the hyperbolic case, it  suffices to prove the analogue of
the results of \cite{CCT2}. We resume the semi-classical approach of
these results, as in \cite[Appendix]{BGTENS} and \cite[Appendix~B]{CaARMA}. 
\smallbreak

We give the proof for a nonlinearity which may be focusing: let
$\kappa\in \{-1,+1\}$. 
Assume $s<s_c$,
and for $0<\l \le 1$, consider:
\begin{equation}
  \label{eq:NLSsc}
  i\d_t u + \Delta_{\H^n} u=\kappa|u|^{2\si}u \quad ;\quad
  u_{\mid t=0}=\l^{-\frac{n}{2}+s}a_0 \(\frac{r}{\l}\)\, , 
\end{equation}
where $a_0 \in C_0^\infty(\R_+;\C)$ is such that 
\begin{equation*}
\operatorname{supp} a_0 \subset \left\{ 1\le r\le 2\right\}.
\end{equation*}
Then $u(0,\cdot)$ is bounded in $H^s(\H^n)$, uniformly for $\l \in ]0,
1]$. Since $u(0,\cdot)$ is radially symmetric, so is $u(t,\cdot)$, and
we write $u(t,r)$. Define $\psi$ by:
\begin{equation*}
  \psi(t,r)= \l^{\frac{n}{2}-s}u\( \l^{\frac{n\si}{2} +1-s\si}t,\l
  r\), \text{ or }u(t,r) = \l^{-\frac{n}{2}+s} \psi\(
  \frac{t}{\l^{\frac{n\si}{2} +1-s\si}},\frac{r}{\l}\).
\end{equation*}
Denote $\eps = \l^{\frac{n\si}{2} -1-s\si}$. Because we assumed $s<s_c$,
$\lambda$ and $\eps$ go to zero simultaneously. We have:
\begin{equation*}
  %\label{eq:psi}
  i\eps \d_t \psi^\eps +\eps^2\d_r^2 \psi^\eps +
  (n-1)\eps^2 \l 
  \coth(\l r) \d_r \psi^\eps =\kappa |\psi^\eps|^{2\si}\psi^\eps\quad ;\quad
  \psi^\eps (0,r)=a_0(r)\, .
\end{equation*}
The idea is that for very small times, the Laplacian is
negligible. Introduce the approximate solution
\begin{equation*}
  %\label{eq:varphi}
  i\eps \d_t \varphi^\eps =\kappa |\varphi^\eps|^{2\si}\varphi^\eps\quad
  ;\quad \varphi^\eps   (0,r)=a_0(r)\, .
\end{equation*}
We have explicitly:
\begin{equation*}
  \varphi^\eps(t,r) = a_0(r)e^{-i \kappa \frac{t}{\eps} |a_0(r)|^2}.
\end{equation*}
In particular, the support of $\varphi^\eps(t,\cdot)$ is the same as that
of $a_0$. 

\begin{prop}\label{prop:OGNL}
   Fix $k>n/2$. Then we can find
  $c_0,c_1,\theta,C>0$ 
  independent of $\eps\in 
  ]0,1]$ such that $\psi^\eps$ and $\varphi^\eps$ satisfy:
\begin{equation*}
 \|\psi^\eps-\varphi^\eps\|_{L^\infty([0,c_0 \eps|\ln
 \eps|^\theta];H^k)}  \le C \eps 
  |\ln \eps|^{c_1}\, . 
\end{equation*}
\end{prop}
\begin{proof}
Denote $w^\eps=\psi^\eps -\varphi^\eps$. It solves:
\begin{align*}
  i\eps \d_t w^\eps +\eps^2\d_r^2 w^\eps + (n-1)\eps^2 \l
  \coth(\l r) \d_r w^\eps =&\kappa \(F(\psi^\eps) -
  F(\varphi^\eps)\) -\eps^2\d_r^2\varphi^\eps \\
&- (n-1)\eps^2 \l
  \coth(\l r) \d_r\varphi^\eps  ,
\end{align*}
 with $ w^\eps_{\mid t=0}=0$, where we have set $F(z)=|z|^{2\si}z$. 
Introduce the vector-fields
\begin{equation*}
  H_j = \om_j \d_r\, .
\end{equation*}
They commute with the Laplacian $\Delta_{\H^n}$. Moreover, since $w$
is radially symmetric,
\begin{align*}
H_j\big(  \d_r^2  + (n-1) &\l
  \coth(\l r) \d_r\big) w^\eps=\\
&= H_j\(  \d_r^2  + (n-1) \l
  \coth(\l r) \d_r+\frac{\l^2}{\sinh^2(\l r)}\Delta_{\Sph^{n-1}}\) w^\eps  \\
&= \(  \d_r^2  + (n-1) \l
  \coth(\l r) \d_r+\frac{\l^2}{\sinh^2(\l r)}\Delta_{\Sph^{n-1}}\) H_j
  w^\eps\, .
\end{align*}
For $k\ge 0$, we apply $H_{j_1}\circ\ldots\circ H_{j_k}$ to the equation
solved by $w^\eps$. The usual $L^2$ estimate yields:
\begin{align*}
  \|w^\eps\|_{L^\infty([0,t];H^k)} \lesssim& \frac{1}{\eps}\left\|F(w^\eps
  +\varphi^\eps) -  F(\varphi^\eps)  \right\|_{L^1([0,t];H^k)}\\
& +\eps
  \|\d_r^2\varphi^\eps\|_{L^1([0,t];H^k)} + \eps \|\l
  \coth(\l r) \d_r\varphi^\eps\|_{L^1([0,t];H^k)}   .
\end{align*}

Since $\operatorname{supp} \varphi^\eps(t,\cdot)\subset \{1\le r\le 2\}$, 
\begin{equation*}
  \|\l
  \coth(\l r) \d_r\varphi^\eps\|_{L^1([0,t];H^k)} \lesssim
  \|\d_r\varphi^\eps\|_{L^1([0,t];H^k)} \, ,
\end{equation*}
and we infer:
\begin{align*}
  \|w^\eps\|_{L^\infty([0,t];H^k)} &\lesssim \frac{1}{\eps}\left\|F(w^\eps
  +\varphi^\eps) -  F(\varphi^\eps)  \right\|_{L^1([0,t];H^k)} +\eps
  \|\varphi^\eps\|_{L^1([0,t];H^{k+2})}\\
& \lesssim \frac{1}{\eps}\left\|F(w^\eps
  +\varphi^\eps) -  F(\varphi^\eps)  \right\|_{L^1([0,t];H^k)} +\eps \int_0^t\<
  \frac{s}{\eps}\>^{k+2}ds. 
\end{align*}
Since $\si$ is an integer, the fundamental theorem of calculus yields, when
$k>n/2$: 
\begin{align*}
 \left\|F(w^\eps(t)
  +\varphi^\eps(t)) -  F(\varphi^\eps(t))  \right\|_{H^k}&\lesssim \(
  \|w^\eps(t)\|^{2\si}_{H^k} +
  \|\varphi^\eps(t)\|^{2\si}_{H^k}\)\|w^\eps(t)\|_{H^k} \\%\label{eq:algebre}\\
&\lesssim \(
  \|w^\eps(t)\|^{2\si}_{H^k} +
  \left\langle\frac{t}{\eps}\right\rangle^{2\si k}\)\|w^\eps(t)\|_{H^k}.%\notag
\end{align*}
On any time interval where we have, say, $\|w^\eps\|_{H^k}\le 1$, we infer:
\begin{equation*}
  \|w^\eps\|_{L^\infty([0,t];H^k)} \le
  \frac{C_0}{\eps}\int_0^t \left\langle\frac{s}{\eps}\right\rangle^{2\si k} 
  \| w^\eps (s)\|_{H^k}ds + C_1\eps \int_0^t \<
  \frac{s}{\eps}\>^{k+2}ds.
\end{equation*}
Gronwall lemma yields:
\begin{align*}
  \|w^\eps\|_{L^\infty([0,t];H^k)} &\lesssim \int_0^{t} \eps  \<
  \frac{s}{\eps}\>^{k+2}\exp\( C_0\int_s^t
  \frac{1}{\eps}\left\langle\frac{\tau}{\eps}\right\rangle^{2\si k}
  d\tau\)ds\\
&  \lesssim \int_0^{t} \eps  \<
  \frac{s}{\eps}\>^{k+2}\exp\( C'
  \frac{t}{\eps}\left\langle\frac{t}{\eps}\right\rangle^{2\si k}\)ds\\
&  \lesssim  \eps t  \<
  \frac{t}{\eps}\>^{k+2}\exp\( C'
  \frac{t}{\eps}\left\langle\frac{t}{\eps}\right\rangle^{2\si k}\).
\end{align*}
For $t = c_0 \eps |\ln \eps|^\theta$, we have:
\begin{align*}
  \|w^\eps\|_{L^\infty([0,c_0 \eps |\ln \eps|^\theta];H^k)} & \lesssim  \eps^2
  |\ln \eps|^{(k+3)\theta}   \exp\( C'c_0^{1+2\si k}
  |\ln \eps|^{(1+2\si k)\theta}\).
\end{align*}
Now if we take $\theta=(1+2\si k)^{-1}$ and $c_0$ sufficiently small,
 we have the estimate of the proposition. A continuity argument
 completes the proof, for $\eps$ sufficiently small.  
\end{proof}
\begin{cor}
  Let $n\ge 2$, $\kappa\in \R\setminus \{0\}$ and $\si>0$.
For $s<\frac{n}{2}-\frac{1}{\si}$, \eqref{eq:nls}
is not locally well-posed in $H^s(\H^n)$: for any
$\delta >0$, we can find families $(u_{01}^\eps)_{0<\eps\le 1}$ and
$(u_{02}^\eps)_{0<\eps\le 1}$ of radially symmetric functions  with
\begin{equation*}
  u_{0j}^\eps(r)\in C_0^\infty(\R_+)\quad ;\quad
  \|u_{01}^\eps\|_{H^s}, 
  \|u_{02}^\eps\|_{H^s} \le \delta\, ,\quad \|u_{01}^\eps -
  u_{02}^\eps\|_{H^s}\to 0 \text{ as }\eps \to 0\, ,
\end{equation*}
such that if $u_1^\eps$ and $u_2^\eps$ denote the solutions to
\eqref{eq:nls} with these initial data, there exist $0<t^\eps\to 0$,
and $c>0$ independent of $\eps\in ]0,1]$, such that 
\begin{equation*}
\left\|u_1^\eps\(t^\eps\)-u_2^\eps\(t^\eps\) \right\|_{H^s}\ge c.
\end{equation*}
\end{cor}
\begin{rema}
  We could also prove the \emph{norm inflation} phenomenon, as called
  in \cite{CCT2}. It is rather this result which is proven in
  \cite[Appendix]{BGTENS}. 
\end{rema}
\begin{proof}
  Let $a_0$ as above, and $u_1$ the solution to \eqref{eq:NLSsc}. Let
  \begin{equation*}
    u_{02}^\eps (r)= (1+\delta^\eps)u_{01}^\eps (r),
  \end{equation*}
where $\delta^\eps \to 0$ as $\eps \to 0$. Denote
\begin{equation*}
  \psi^\eps_j(t,r)= \l^{\frac{n}{2}-s}u_j\( \l^{\frac{n\si}{2} +1-s\si}t,\l
  r\).
\end{equation*}
>From Proposition~\ref{prop:OGNL}, for $k>n/2$, 
\begin{equation*}
  \left\| \psi^\eps_j - \varphi^\eps_j\right\|_{L^\infty([0,c_0 \eps |\ln
  \eps|^\theta];H^k)} \to 0 \text{ as }\eps\to 0\, .
\end{equation*}
Note that the constants $c_0$ and $\theta$ for $\psi_2$ can be taken
uniform with respect to $\eps$, since $\delta^\eps$ is bounded.
Now we have
\begin{align*}
  \left\| \varphi^\eps_1(t) - \varphi^\eps_2(t)\right\|_{\dot H^s} &=\left\|
  a_0 e^{i\kappa \frac{t}{\eps}|a_0|^2} - \(1+\delta^\eps\) a_0 e^{i\kappa
  \frac{t}{\eps}\(1+\delta^\eps\)^2|a_0|^2}\right\|_{\dot H^s}   \\
&\Eq \eps 0 \left\|
  a_0 \( e^{i\kappa \frac{t}{\eps}|a_0|^2} - e^{i\kappa
  \frac{t}{\eps}\(1+\delta^\eps\)^2|a_0|^2}\)\right\|_{\dot H^s}\\
&\Eq \eps 0 \left\|
  a_0 e^{i\kappa \frac{t}{\eps}|a_0|^2}\( 1 - e^{i\kappa
  \frac{t}{\eps}\(2\delta^\eps+(\delta^\eps)^2\)|a_0|^2}\)\right\|_{\dot H^s}.
\end{align*}
With $t^\eps = c_0 \eps |\ln \eps|^\theta$, we find:
\begin{align*}
  \left\| \varphi^\eps_1(t^\eps) -
  \varphi^\eps_2(t^\eps)\right\|_{\dot H^s} \Eq \eps
  0 \left\| 
  a_0 \( 1 - e^{i\kappa
  2\delta^\eps |\ln \eps|^\theta |a_0|^2}\)\right\|_{\dot H^s}.
\end{align*}
If we take $\delta^\eps =|\ln \eps|^{-\theta}$, then the corollary
follows, since $\psi$ and $u$ have the same $\dot H^s$  norms. 
\end{proof}

\section{Asymptotic behavior of free solutions in $L^2_{\rm rad}(\H^3)$}
\label{sec:DA}
We stick to the case $n=3$, because the Harish-Chandra coefficient is
simpler. In this appendix, all irrelevant
``physical''/geometrical constants are denoted by $c$. 
\smallbreak

The Fourier transform in the radially symmetric case is defined as:
\begin{equation}
  \label{eq:Fourier3D}
  \widehat f(\l) = \frac{c}{\l}\int_0^\infty \sin(\l r)f(r)\sinh r
  dr\, .
\end{equation}
Plancherel formula reads:
\begin{equation}
  \label{eq:Plancherel3D}
  \int_0^\infty |f(r)|^2 \sinh^2 r dr = c \int_0^\infty \left|
  \widehat f(\l)\right|^2 \l^2 d\l\, . 
\end{equation}
\begin{lem}\label{lem:DA}
  For $u_0\in L^2 (\H^3)$, radially symmetric, denote:
  \begin{equation*}
    u_{\rm asym}(t,r) = c\frac{e^{-it
    +i\frac{r^2}{4t}}}{t^{3/2}}\frac{r}{\sinh r} \widehat
    {u_0}\(\frac{r}{2t}\).
  \end{equation*}
Then we have:
  \begin{equation*}
    \left\| e^{it\Delta_{\H^3}}u_0 - u_{\rm
    asym}(t)\right\|_{L^2 (\H^3)}\Tend t {+\infty} 0 \, . 
  \end{equation*}
\end{lem}
\begin{proof}
  First, we show an explicit representation for radial solutions:
  \begin{equation}
  \label{eq:solradial}
   e^{it\Delta_{\H^3}}u_0(r) =  c\frac{e^{-it
    +i\frac{r^2}{4t}}}{t^{3/2}} \int_0^\infty
    \frac{e^{i\frac{\rho^2}{4t}}}{\sinh r \sinh \rho} t\sin \(
    \frac{r\rho}{2 t}\) u_0(\rho) \sinh^2 \rho d \rho \, .
  \end{equation}
To prove this we recall the representation of the free solution for $n=3$:
\begin{equation}
  \label{eq:n=3}
  u(t,\Omega)= \frac{c}{|t|^{3/2}} e^{-it}\int_{\H^3}
  u_0(\Omega')\,e^{i\frac{d^2(\Omega,\Omega')}{4t}}
\frac{d(\Omega,\Omega')}{\sinh d(\Omega,\Omega')}\,d\Omega'\,, 
\end{equation}
where $d(\Omega,\Omega')$ is the hyperbolic distance between $\Omega$ and
$\Omega'$ (see \cite{BaHyper}). 
>From \eqref{eq:n=3}, one gets that for radial initial data, the free
solution writes 
 \begin{equation*}
   e^{it\Delta_{\H^3}}u_0(\cosh r, \sinh r\omega) =
   c\frac{e^{-it}}{t^{3/2}} \int_0^\infty    
   \int_{\Sph^2} K(t,r,\rho,\omega\cdot\omega') d\omega' u_0(\rho) \sinh^2
   \rho d \rho \, , 
  \end{equation*}
  with 
  \begin{equation*}
  K(t,r,\rho,\omega\cdot\omega')=e^{i\frac{z^2}{4t}}
  \frac{z}{\sinh z} \Big|_{z=\cosh^{-1}\,(\cosh r\cosh\rho-\sinh
  r\sinh\rho\,\omega\cdot\omega')}\,. 
  \end{equation*}
  Let us consider an isometry $T\in SO(3)$ such that
  $T(1,0,0)=\omega$. Then a given $\omega'\in\R^3\cap \Sph^2$ defines a
  unique pair $(\alpha,\theta)\in(0,\pi)\times\R^2\cap \Sph^1$,
  related by the formula:
   \begin{equation*}
  \omega'=T(\cos\alpha,\sin\alpha\,\theta).
  \end{equation*}
  Moreover, $\omega\cdot\omega'=\cos\alpha$ and $d\omega'=\sin\alpha\,
  d\alpha\, d\theta$. With this change of variable, 
  \begin{equation*}
  \int_{\Sph^2} K(t,r,\rho,\omega\cdot\omega') d\omega'=c\int_0^\pi
  K(t,r,\rho,\cos\alpha) \sin\alpha\,d\alpha \,. 
  \end{equation*}
  Next, we change $\cos\alpha$ into $x$, so
  \begin{equation*}
  \int_{\Sph^2} K(t,r,\rho,\omega\cdot\omega') d\omega'=c\int_{-1}^1
  K(t,r,\rho,x)\,dx \,. 
  \end{equation*}
  Finally, we do a last change of variable, 
  \begin{equation*}
  \cosh r\cosh\rho-\sinh r\sinh\rho\,x=\cosh y\, ,  
  \end{equation*}
  and get
  \begin{equation*}
  \int_{\Sph^2} K(t,r,\rho,\omega\cdot\omega') d\omega'=\frac{c}{\sinh
  r\,\sinh\rho}\int_{|r-\rho|}^{r+\rho} e^{i\frac{y^2}{2t}} y \,dy. 
  \end{equation*}
  By simple integration formula (\ref{eq:solradial}) follows. We infer:
\begin{align*}
  e^{it\Delta_{\H^3}}u_0(r) - &u_{\rm
    asym}(t,r)=\\
&= c\frac{e^{-it
    +i\frac{r^2}{4t}}}{t^{3/2}} \frac{r}{\sinh r} \int_0^\infty
    \(e^{i\frac{\rho^2}{4t}}-1\) \frac{2t}{r}\sin \(
    \frac{r\rho}{2t}\) u_0(\rho) \sinh \rho d \rho\, .
\end{align*}
Therefore:
\begin{align*}
  \Big\|e^{it\Delta_{\H^3}}u_0 - &u_{\rm
    asym}(t)\Big\|^2_{L^2}=\\
&= \frac{c}{t^{3}} \int_0^\infty r^2 \left| \int_0^\infty
    \(e^{i\frac{\rho^2}{4t}}-1\) \frac{t}{r}\sin \(
    \frac{r\rho}{2t}\) u_0(\rho) \sinh \rho d \rho\right|^2 dr\\
& = c \int_0^\infty r^2 \left| \int_0^\infty
    \(e^{i\frac{\rho^2}{4t}}-1\) \frac{1}{r}\sin \(r\rho\) u_0(\rho)
    \sinh \rho d \rho\right|^2 dr \, .
\end{align*}
Use Plancherel formula \eqref{eq:Plancherel3D}:
\begin{equation*}
  \Big\|e^{it\Delta_{\H^3}}u_0 - u_{\rm
    asym}(t)\Big\|^2_{L^2}
 = c   \int_0^\infty
    \left|\(e^{i\frac{\rho^2}{4t}}-1\)  u_0(\rho)\right|^2
    \sinh^2 \rho d \rho  \, .
\end{equation*}
Now we conclude with a density argument, thanks to the estimate (for
instance) 
\begin{equation*}
  |e^{i\theta} -1|\lesssim |\theta|\, . 
\end{equation*}
\end{proof}
\begin{rema}
  This asymptotic behavior is essentially the same as in the
  Euclidean case, up to a new
  oscillation in time, and the weight 
  $\frac{r}{\sinh r}$. This can also be seen as follows: in the proof,
  we have used the identity
  \begin{equation*}
   e^{it\Delta_{\H^3}}u_0(r) = W \M_t D_t \F \M_t u_0\, , 
  \end{equation*}
where  $\F$ is the Fourier transform, $\M_t(r)$ is the multiplication
by $e^{ir^2/(4t)}$, $D_t$ is the  dilation at scale $1/(2t)$ with
$L^2$ scaling:
\begin{equation*}
  D_t\varphi(r) = \frac{1}{t^{3/2}}\varphi\(\frac{r}{2t}\),
\end{equation*}
and $W$ is the weight $ W= e^{-it}\frac{r}{\sinh r}$.
In the Euclidean case, we have the same formula, with the only change
$W={\rm Id}$, and the usual asymptotics is:
\begin{equation*}
  e^{it\Delta_{\R^n}} = \M_t D_t \F \M_t\Eq t {+\infty}\M_t
  D_t \F\, .
\end{equation*}
The proof of the $L^2$ asymptotics is the same as above. 
\end{rema}

\section{A Galilean operator?}
\label{sec:galileo}

In the Euclidean case, a nice object for scattering theory (and also
blow-up, see e.g. \cite{Weinstein83}) is the Galilean operator
\begin{equation*}
  J_{\rm eucl}(t) =x +2it \nabla_x = 2it e^{i\frac{|x|^2}{4t}}\nabla_x
  \(e^{-i\frac{|x|^2}{4t}} \cdot \)= e^{it\Delta_{\R^n}} x
  e^{-it\Delta_{\R^n}} \, . 
\end{equation*}
Recall that for any \emph{radial} function $\phi$ on $\R^3$, we have:
\begin{equation*}
  \Delta_{\R^3} \(\frac{\phi(r)}{r}\) = \frac{1}{r}\d_r^2\phi. 
\end{equation*}
A similar identity is available for radial functions on $\H^3$:
\begin{equation*}
  \Delta_{\H^3} \(\frac{\phi(r)}{\sinh r}\) = \frac{1}{\sinh r}\d_r^2\phi. 
\end{equation*}
Using the Galilean operator on the (half-)line, it is then natural to
introduce the following operator, acting on radial functions on
$\H^3$:  
\begin{equation}\label{eq:J}
  J(t) =  \frac{2it}{\sinh r} e^{i\frac{r^2}{4t}}\d_r \(
  e^{-i\frac{r^2}{4t}} \sinh r \, \cdot\)=r+2it\d_r +2it
  \coth r \, .
\end{equation}
Now we have
\begin{equation*}
  \left[J,i\d_t
+\Delta_{\H^3}\right]u(t,r)=0\, ,
\end{equation*}
provided that $u$ is radial. Note that as in the Euclidean
case, $J$ is an Heisenberg observable (see e.g. \cite{Robert}): 
\begin{equation}\label{eq:Heisenberg}
  J(t) = e^{it\Delta_{\H^3}} r
  e^{-it\Delta_{\H^3}}\, .
\end{equation}
We already saw that in the radial framework, $J$ commutes with the
Schr\"odinger operator. Proceeding as in \cite[Lemma~6.2]{CM}, we find:
\begin{lem}\label{lem:dispJ}
  Let $n=3$. For every $q\in [2,6[$, there exists $c_q$ such that
  \begin{equation*}
    \left\|{\tt w}_3^{1-\frac{2}{q}}\phi \right\|_{L^q}\le
  \frac{c_q}{|t|^{\delta(q)}} \|\phi\|_{L^2}^{1-\delta(q)}
  \|J(t)\phi\|_{L^2}^{\delta(q)},\ \forall t\not =0,\quad \text{for every
  \emph{radial} function }\phi,   
  \end{equation*}
where we have denoted 
$\displaystyle {\tt w}_3(r)=\frac{\sinh r}{r}$ and $\displaystyle \delta(q) =
 3\(\frac{1}{2}-\frac{1}{q}\)\in [0,1]$. 
\end{lem}
It is important to understand how $J$ acts on nonlinear terms. Let $F$
be a $C^1$ function such that $F(z)=G(|z|^2)z$ (the usual gauge
invariance). We compute:
\begin{equation}
  \label{eq:JNL}
  J(t)F(u) = \d_z F(u)J(t)u - \d_{\bar z} F(u)\overline{J(t)u} +it
  \coth r F(u). 
\end{equation}
Forgetting the last term, we would have the same expression  as in $\R^n$,
and $J$ would act on such nonlinearities like a
derivative. Unfortunately, this last term cumulates two features:
extra linear growth in time, and singularity as $r\to 0$. 
\smallbreak

As a matter of fact, the above drawback is also present in the
radial Euclidean case. There, it can be removed by using $J_{\rm
  eucl}$, even in a radial framework (see
e.g. \cite{BourgainRadial}). However, the analogue for $J_{\rm 
  eucl}$ in hyperbolic space (not only in the radial case) may
just not exist\ldots
Consider the Euclidean case in $\R^3$. We have 
\begin{equation*}
  \Delta_{\R^3}= \d_r^2 + \frac{2}{r}\, \d_r
  +\frac{1}{r^2}\Delta_{\Sph^{2}}\, . 
\end{equation*}
For a radial function $u(t,r)$, we have the commutation relation
\begin{equation*}
  \left[x_j+2it\d_j, i\d_t+\Delta_{\R^3} \right]u(t,r)=0,
\end{equation*}
where 
$  (x_j+2it\d_j)u(t,r)= (r\om_j +2it\om_j \d_r)u(t,r).$
The factor $\om_j$ is crucial: in general,
\begin{equation*}
  \left[r+2it\d_r, i\d_t+\Delta_{\R^3} \right]u(t,r)\not=0.
\end{equation*}
On the other hand,
\begin{equation*}
  \left[r+2it\d_r+2i\frac{t}{r}, i\d_t+\Delta_{\R^3} \right]=0.
\end{equation*}
The above operator can also be written as 
\begin{equation*}
  {\bf J}(t)= r+2it\d_r+2i\frac{t}{r} =  \frac{2it}{r}
  e^{i\frac{r^2}{4t}}\d_r \( 
  e^{-i\frac{r^2}{4t}} r \, \cdot\)\, .
\end{equation*}
When acting on gauge invariant nonlinearities, it is not like a
derivative:
\begin{equation*}
  {\bf J}(t)F(u) = \d_z F(u){\bf J}(t)u - \d_{\bar z}
  F(u)\overline{{\bf J}(t)u} +\frac{it}{r} F(u). 
\end{equation*}
The last factor has the same drawback as above. Note that even in the
radial case, one uses $J_{\rm eucl}$ and not ${\bf J}$ (see
e.g. \cite{BourgainRadial}). So we may try to extend $J$ to a
non-radial framework. Yet, seeking $J_{\rm hyper}$ of the form
\begin{equation*}
  J_{\rm hyper} = r\om_j +2it\om_j \d_r +2it h(r) =
  \frac{2it}{f(r)}e^{i\frac{r^2}{4t}}\om_j\d_r \( 
  e^{-i\frac{r^2}{4t}} f(r) \, \cdot\)\quad \(h=\frac{f'}{f}\),
\end{equation*}
and writing 
\begin{equation*}
  \left[J_{\rm hyper}, i\d_t+\Delta_{\H^3} \right]u(t,r)=0,
\end{equation*}
yields incompatible conditions. 
\smallbreak 

Using the operator $J$ and Lemma~\ref{lem:dispJ}, one can prove the
following scattering result 
though. For $T>0$ possibly large, define:
\begin{equation*}
\begin{aligned}
  Y_T=\{& u,Ju\in C(]-\infty,-T];L^2),\text{ with }
  {\tt w}_3^{2/3}u,{\tt w}_3^{2/3}Ju \in L^2(]-\infty,-T];L^6)\ ;\\
 & \|u\|_{L^\infty(]-\infty,-T;L^2)} \le 2 \|u_-\|_{L^2},\  
 \|{\tt w}_3^{2/3}u\|_{L^2(\R;L^6)} \le 2 C_6\|u_-\|_{L^2} ,\\
& \|Ju\|_{L^\infty(]-\infty,-T;L^2)} \le 2 \|ru_-\|_{L^2},\  
 \|{\tt w}_3^{2/3}Ju\|_{L^2(\R;L^6)} \le 2 C_6\|ru_-\|_{L^2} \},
\end{aligned}
 \end{equation*}
where $C_6$ is given by Proposition~\ref{prop:weightedstri3+} when
$n=3$. 
\begin{prop}\label{prop:waveJ}
 Let $n=3$, $t_0=-\infty$ and $1/4<\si<1$. For every
  \emph{radial} function $\varphi=u_-\in L^2(\H^3)$ with $r u_-\in
  L^2(\H^3)$, there 
  exists $T=T(\si,  \|u_-\|_{L^2},\|ru_-\|_{L^2})$ such that
  \eqref{eq:duhamel} has a unique solution in $Y_T$. 
\end{prop} 
\begin{rema}
  The condition $1/4<\si<1$ looks rather strange at first glance. It
  appears because of the singular term that shows up when $J$ acts on
  the nonlinearity, as discussed above.  Without this term, we could
  virtually cover the 
  range $0<\si <2$. Note however that we can go below $\si=1/3$, thus
  showing the absence of (the usual) long range effects.
\end{rema}
\begin{proof}
  We want to show that the map 
  \begin{equation*}
    \Phi(u)(t):= U(t)u_- -i \int_{-\infty}^t U(t-s)|u|^{2\si}u(s)ds
  \end{equation*}
has a fixed point in $Y_T$ for $T$ sufficiently large. Let $(q,p)$ be
an $3$-admissible pair to be chosen
later. Proposition~\ref{prop:weightedstri3+} yields:
\begin{equation*}
  \left\| {\tt w}_3^{2/3}\Phi(u)\right\|_{L_T^2L^6}\le C_6
  \|u_-\|_{L^2} + C_{6,p}\left\|
  {\tt w}_3^{1-2/{p'}}|u|^{2\si}u\right\|_{L_T^{q'}L^{p'}},
\end{equation*}
where we denote from now on: $L^a_TL^b :=
L^a(]-\infty,-T];L^b(\H^3))$. 
Introduce indices such that:
\begin{equation}\label{eq:holderJ}
\left\{
  \begin{aligned}
    \frac{1}{q'} &= \frac{1}{q} +\frac{2\si}{s}+\frac{1}{\theta}\virgp\\
    \frac{1}{p'}&=  \frac{1}{p} + \frac{2\si}{k}\virgp
  \end{aligned}
\right. \qquad \text{ with }s\in ]2,6[.
\end{equation}
H\"older's inequality then yields:
\begin{equation*}
 \left\|
  {\tt w}_3^{1-2/q'}|u|^{2\si}u\right\|_{L^{p'}_TL^{q'}} \le  \left\|
  {\tt w}_3^{1-2/q}u\right\|_{L^p_TL^{q}} \left\|
  {\tt w}_3^{1-2/s}u\right\|_{L^k_TL^s}^{2\si}\left\| {\tt w}_3^{4/q-2 
  +2\si(2/s-1)}\right\|_{L^\theta}. 
\end{equation*}
The last term is finite provided that:
\begin{equation}\label{eq:thetaJ}
  2-\frac{4}{q} +2\si \( 1-\frac{2}{s}\)
  >\frac{2}{\theta}\cdot 
\end{equation}
In view of \eqref{eq:holderJ}, this is equivalent to $\si>0$.

Now we show that \eqref{eq:holderJ} can be achieved
for a  $3$-admissible pair $(p,q)$. Letting $(p,q)=(2,6)$, we find
$k=\infty$ and \eqref{eq:thetaJ} is satisfied for any $\si>0$. Then
$s\in ]2,6[$ provided that  $\si <2$. Since these algebraic conditions
are open, they still hold if we take $q = 6-\eps$ for $\eps>0$
sufficiently small, and $p$ such that $(p,q)$ is $3$-admissible:
 \begin{equation*}
 \left\|
  {\tt w}_3^{1-2/q'}|u|^{2\si}u\right\|_{L^{p'}_TL^{q'}} \le  C \left\|
  {\tt w}_3^{1-2/q}u\right\|_{L^p_TL^{q}} \left\|
  {\tt w}_3^{1-2/s}u\right\|_{L^k_TL^s}^{2\si}. 
\end{equation*}
Note also that now, $p>2$, hence $k<\infty$. 
By interpolation, the first factor of the right hand side is
controlled by:
\begin{equation*}
  \left\|
  {\tt w}_3^{1-2/q}u\right\|_{L^p_TL^{q}} \le \left\|
  u\right\|_{L^\infty_TL^2}^{1-\delta(q)}\left\|
  {\tt w}_3^{2/3}u\right\|_{L^2_TL^{6}}^{\delta(q)}. 
\end{equation*}
>From Lemma~\ref{lem:dispJ}, we have, for $t\le -T$:
\begin{equation*}
  \left\|
  {\tt w}_3^{1-2/s}u(t)\right\|_{L^s}\lesssim
  \frac{1}{|t|^{\delta(s)}}\|u\|_{L^\infty_TL^2}^{1-\delta(s)}
\|Ju\|_{L^\infty_TL^2}^{\delta(s)}.  
\end{equation*}
We check that \eqref{eq:holderJ} implies
\begin{align*}
  k\delta(s) &= 2 +\frac{k}{\si}\( -1 +\frac{3\si}{2}
  +\frac{3}{2}\(1-\frac{2}{q}\) -\frac{3\si}{s}\)\\
& > 2 + \frac{k}{\si}\( \frac{3\si}{2}
   -\frac{3\si}{s}\)>2. 
\end{align*}
Therefore, 
\begin{equation*}
  \left\|
  {\tt w}_3^{1-2/s}u\right\|_{L^k_TL^s}^{2\si} \lesssim T^{-2\si
  /k}\(\|u\|_{L^\infty_TL^2} +  \|Ju\|_{L^\infty_TL^2}\)^{2\si}. 
\end{equation*}
Then choosing $T$ sufficiently large, we see that for $u\in Y_T$, 
\begin{equation*}
 \left\| {\tt w}_3^{2/3}\Phi(u)\right\|_{L_T^2L^6}\le2 C_6
  \|u_-\|_{L^2}. 
\end{equation*}
The similar estimate for $\| \Phi(u)\|_{L_T^\infty L^2}$
proceeds along the same lines.
\smallbreak

To estimate $J \Phi(u)$, we find:
\begin{align*}
  \left\| {\tt w}_3^{2/3}J\Phi(u)\right\|_{L_T^2L^6}\le& \ C_6
  \|r u_-\|_{L^2} + C\left\|
  {\tt w}_3^{1-2/q'}|u|^{2\si}Ju\right\|_{L_T^{p'}L^{q'}}\\
& + C\left\|
  {\tt w}_3^{1-2/q_1'}t\coth r |u|^{2\si+1}\right\|_{L_T^{p_1'}L^{q_1'}},
\end{align*}
with the same admissible pair $(p,q)$ as before, and where $(p_1,q_1)$
is a possibly different admissible pair. For the second term
of the right hand side, we proceed as before, to find:
 \begin{equation*}
 \left\|
  {\tt w}_3^{1-2/q'}|u|^{2\si}Ju\right\|_{L^{p'}_TL^{q'}} \le  C \left\|
  {\tt w}_3^{1-2/q}Ju\right\|_{L^p_TL^{q}} \left\|
  {\tt w}_3^{1-2/s}u\right\|_{L^k_TL^s}^{2\si}. 
\end{equation*}
We are left with the next term, involving ${\tt w}_3^{1-2/q'}t\coth r
|u|^{2\si+1}$. Introduce the condition
\begin{equation}\label{eq:holdersing}
    \frac{1}{q_1'} = \frac{1}{q_1}
    +\frac{2\si}{s_1}+\frac{1}{\theta_1}\cdot
\end{equation}
When it is satisfied for $s_1,\theta_1\ge 1$, we have:
\begin{align*}
  \left\|
  {\tt w}_3^{1-2/q_1'}\coth r |u|^{2\si+1}\right\|_{L^{q_1'}}\le &\, \left\|
  {\tt w}_3^{1-2/q_1}u\right\|_{L^{q_1}} \left\|
  {\tt w}_3^{1-2/s_1} u\right\|_{L^{s_1}}^{2\si}\times\\
&\times\left\|
  {\tt w}_3^{4/q_1 -2 +2\si(2/s_1-1)}\coth r \right\|_{L^{\theta_1}}.
\end{align*}
For the last term to be finite, a new condition appears, for the
integral to converge near $r=0$ (for $r\to \infty$, nothing is
changed): 
\begin{equation}\label{eq:theta1}
  1\le\theta_1<3. 
\end{equation}
For $q_1,s_1\in ]2,6[$, Lemma~\ref{lem:dispJ} then yields:
\begin{align*}
  \left\|
  {\tt w}_3^{1-2/q_1'}\coth r |u|^{2\si+1}\right\|_{L^{q_1'}}\lesssim
  \frac{1}{|t|^{\delta(q_1)+ 2\si \delta(s_1)}}\(
  \|u\|_{L^\infty_TL^2} +  \|Ju\|_{L^\infty_TL^2}\)^{2\si+1}. 
\end{align*}
The right hand side multiplied by $t$ is in $L^{p_1'}(]-\infty,-T])$
as soon as: 
\begin{equation}\label{eq:int1}
  \delta(q_1)+ 2\si \delta(s_1)-1>\frac{1}{p_1'}\cdot
\end{equation}
So we are left with the following situation: if we can meet
\eqref{eq:holdersing}, \eqref{eq:theta1} and \eqref{eq:int1} with
$(p_1,q_1)$ admissible and $q_1,s_1\in ]2,6[$, then choosing $T$
sufficiently large, $\Phi$ maps $Y_T$ to itself. 
\smallbreak

The line of reasoning is the same as above: if we pick $q_1=s_1=6$,
then \eqref{eq:holdersing} and \eqref{eq:theta1} imply $\si <1$, and
\eqref{eq:int1} yields $\si >1/4$. Conversely, if $1/4<\si<1$, then
taking $(p_1,q_1)=(2,6)$ and $s_1=6$, $\theta_1$ given by
\eqref{eq:holdersing} satisfies \eqref{eq:theta1}, and
\eqref{eq:int1} holds. By continuity, all the conditions required are
satisfied if we take $q_1 =6-\eps_1$, with $\eps_1>0$ sufficiently small.
\smallbreak

Up to increasing $T$, $\Phi$ is a contraction on  $Y_T$, and the
Proposition~\ref{prop:waveJ} follows. 
\end{proof}
We finally notice that an analogue to the pseudo-conformal
conservation law \cite{GV79Scatt} is available:
\begin{equation}
  \label{eq:Ju}
\begin{aligned}
  \frac{d}{dt}&\(\|J(t)u\|_{L^2}^2 + \frac{4t^2}{\si
  +1}\|u(t)\|_{L^{2\si+2}}^{2\si+2} \)=\\
 =\frac{4t}{\si +1}&\int_0^\infty
  \(2-\si-2\si r\coth r\)|u|^{2\si+2}(t,r)\sinh^2 r dr .
\end{aligned}
\end{equation}
This evolution law should make it possible to establish some
asymptotic completeness results in weighted Sobolev spaces, but we
leave out the discussion here.

\providecommand{\bysame}{\leavevmode\hbox to3em{\hrulefill}\thinspace}

\end{document}